\documentclass[12pt]{article}

\usepackage{latexsym,amsmath,amscd,amssymb,framed,graphics,color,mathrsfs,stmaryrd}
\usepackage{enumerate}

\usepackage{graphicx}

\usepackage[square,authoryear]{natbib}
\definecolor{sand}{rgb}{0.76, 0.7, 0.5}
\definecolor{taupegray}{rgb}{0.55, 0.52, 0.54}
\usepackage[colorlinks]{hyperref}
\hypersetup{%
pdftitle={Field theoretic variational integrator for nonsmooth Lagrangian mechanics},
pdfauthor={Demoures et al.},
pdfsubject={Field theoretic nonsmooth mechanics, variational integrator, contacts.},
pdfstartview={FitH},
urlcolor=sand,
citecolor=blue,
linkcolor=taupegray,
}
\usepackage{url}
\usepackage{framed}
\usepackage{float}
\usepackage{ulem}

\usepackage[all]{xy}

\usepackage[bbgreekl]{mathbbol}

\makeatletter

\@addtoreset{figure}{section}
\def\thefigure{\thesection.\@arabic\c@figure}
\def\fps@figure{h, t}
\@addtoreset{table}{bsection}
\def\thetable{\thesection.\@arabic\c@table}
\def\fps@table{h, t}
\@addtoreset{equation}{section}

\makeatother

\begin{document}

\newtheorem{theorem}{Theorem}[section]
\newtheorem{definition}[theorem]{Definition}
\newtheorem{lemma}[theorem]{Lemma}
\newtheorem{remark}[theorem]{Remark}
\newtheorem{proposition}[theorem]{Proposition}
\newtheorem{corollary}[theorem]{Corollary}
\newtheorem{example}[theorem]{Example}

\def\below#1#2{\mathrel{\mathop{#1}\limits_{#2}}}

\title{
Thermoplasticity \\
as a nonsmooth phenomenon}

\author{Fran\c{c}ois Demoures$^{1,2}$}
\addtocounter{footnote}{1}  \footnotetext{EPFL, Doc \& Postdoc Alumni}
\addtocounter{footnote}{1}\footnotetext{DFJC/DGEP \texttt{francois.demoures@vd.educanet2.ch}}

%\date{ June 2016 }

\maketitle

\makeatother

%\maketitle

%|||-------------------text width----------------------|||
\abstract

This paper develops the variational multisymplectic formulation of nonsmooth elastoplastic phenomena, where the rate of change of plastic strain and the associated thermodynamic entropy evolve by jumps. The formulation relies on convex analysis to describe the plastic non smoothness.

\medskip

\tableofcontents

%%%%%%%%%%%%%%%%%%%%%%%%%%%%%%%%%%%%%%%%%%%%%%%
%%%%%%%%%%%%%%%%%%%%%%%%%%%%%%%%%%%%%%%%%%%%%%%%

\section{Introduction}

\paragraph{Plasticity.} The first important results concerning plasticity are due to \cite{Tresca1872} and \cite{StVenant1871a, StVenant1871b, StVenant1871c}. See \cite{Maugin2016} for an comprehensive historical review, and \cite{Lubliner1990} for a general overview of the subject

The plasticity theory, we consider in this article, was defined in \cite{Hill1950} for the maximal dissipation principle, \cite{GrNa1965} for the formal additive decomposition of the finite Lagrangian strain tensor, \cite{Rice1970, Rice1971} for the theoretical foundations of inelastic constitutive laws for solids, \cite{Rockafellar1970}, \cite{Moreau1973, Moreau1976} for the convex analysis formulation, and \cite{Suquet1979} where the existence of perfect plastic solutions is investigated.   
The classical multiplicative decomposition $\mathbf{F}= \mathbf{F}_e\mathbf{F}_p$ is due to \cite{BiGaSt1957}, \cite{Kroner1960}, \cite{LeLi1967}.

In \cite{SiHu1998}, and \cite{Simo1998} was developed an overview of numerical analysis dedicated to the simulation of problems involving plastic deformation.
More recent complements can be found in, e.g. \cite{Gurtin2000}, \cite{Armero2008}, \cite{ClBa2009}.

Regarding rheological thermodynamics, we refer to the following books and papers: \cite{LaLi1959}, \cite{TrNo1965}, \cite{Truesdell1968}, \cite{ZiWe1987}, \cite{Ottinger2005}, \cite{GuFrAn2010} and \cite{Maugin2011}.
Concerning thermoplasticity, we refer to: \cite{Eckart1948}, \cite{Naghdi1960}, \cite{Ziegler1963}, \cite{GrNa1966}, \cite{Maugin1992}. For numerical study and simulation of thermoplasticity see, for example, \cite{Simo1998} \S $5$. 

From the $18$th century, variational formulations were proposed in order to accommodate irreversible transformations, due to viscosity, plasticity, or heat dissipation, in addition to reversible elastic transformations. The principles involved in these attempts are numerous. We can cite among others Onsager's principle, principle of minimum rate of entropy production, Lagrange-d'Alembert principle, Reissner, Hellinger principles, and Hu-Washizu principles. 
 
Concerning variational principles dedicated to viscoelasticity see \cite{Biot1955}, \cite{FrGe1958}, \cite{Onat1962}, \cite{Gurtin1963}.
Variational principles for plasticity see \cite{Washizu1955, Washizu1968}. 
  Variational principles for thermodynamics see \cite{MaSe2006}, \cite{GBYo2016a, GBYo2016b}. Variational inequalities, see e.g. \cite{EkTe1974}, \cite{GlTa1989}, \cite{HaRe1999}.

However, the plastic domain is defined as the boundary of the elasticity domain, which are jointly formulated by means of nonsmooth inequality constraints. The corresponding dissipated power verifies the maximum-dissipation principle which is equivalent to minimize minus the dissipated power under inequality constraints, i.e., we need to solve an optimization problem, through variational geometry of convex sets. In this vein we can cite \cite{Kachanov1942}, \cite{HoPr1949}, \cite{BuPe1956}, \cite{RoWe1998}.

Concerning the dislocation mechanics involved in crystal plasticity phenomenon, we refer to \cite{MMAVP2008}, \cite{ClMDBa2006}, \cite{FrTaCa2011}, \cite{YaGo2012}.

\paragraph{Multisymplectic formulation.} The plasticity phenomenon will be formulated within the context of multisymplectic continuum mechanics (\cite{GoIsMa2004}) and in particular of multisymplectic nonsmooth continuum mechanics (\cite{FeMaWe2003}). 

\paragraph{Rheological model.} The elastoplastic material exhibits both plastic and elastic behaviour. One can build up a model of nonsmooth elastoplasticity by combining a linear elastic spring and a non-smooth frictional pad. These are known as rheological models originating from the work of \cite{Zener1948}.
See, e.g. the rheological models described in \cite{Maugin1992}, \cite{GuSc1995}, \cite{SiHu1998}, and in \cite{Lion2000}.

\paragraph{Goals and general framework.}
In this paper, firstly we develop the multisymplectic formulation of nonsmooth elastoplastic phenomena in \S \ref{nonsmooth_elastoplasticity}, with a multiplicative decomposition of the total deformation gradient $\mathbf{F}=\mathbf{F}_e\mathbf{F}_p$ into an elastic deformation part $\mathbf{F}_e$ and a plastic deformation part $\mathbf{F}_p$.

Secondly we develop a rheological model dedicated to crystal elastoplasticity and thermoplasticity (with temperature and entropy variables added). 
According to \cite{Simo1998} we admit an additive decomposition of the total strain $\boldsymbol{\epsilon}$ of the system into an elastic strain $\boldsymbol{\epsilon}_e$ and a plastic strain $\boldsymbol{\epsilon}_p$ due to sliding, i.e., $\boldsymbol{\epsilon}=\boldsymbol{\epsilon}_e+\boldsymbol{\epsilon}_p.$
This additive decomposition is consistent with our rheological model composed of a spring and a frictional pad where deformations are small (see Fig. \ref{elast_bar_0}). 

% see Lubliner.

%But as soon as we want to study a realistic problem where is involved plastic materials, all the results of our work can be introduced into the mechanical model to study. Which is generally described by deformation field $\varphi$, with gradient deformation $\mathbf{F}$ and multiplicative decomposition $\mathbf{F}= \mathbf{F}^e\mathbf{F}^p$. To achieve this goal, the starting point for such a development could be \cite{Simo1998}.

\medskip

\noindent \textbf{Summary of the main results:} \\
\noindent $\bullet$ In \S \ref{non_smooth_elastoplastic} the elastoplastic phenomena are expressed through the convex analysis formulation (see \cite{Moreau1973}) and the variational multisymplectic formulation of nonsmooth continuum mechanics developed in \cite{FeMaWe2003}. Where these phenomena are characterized by plastic dissipation jumps. \\
\noindent $\bullet$ In \S\ref{perfect_plastic}, \ref{isotropic_hard}, \ref{kinematic_hard} using the forementioned rheological model, we deduce from the previous results a variational multisymplectic nonsmooth formulation of elastoplasticity with the different following situations: perfect plasticity, isotropic hardening, and kinematic hardening. \\
\noindent $\bullet$ In \S \ref{thermo_perf_plast}, \ref{thermo_isot_hard}, \ref{thermo_kin_hard} assuming perfect plasticity, isotropic hardening, kinematic hardening, and isothermal plasticity, we combine a variational multisymplectic nonsmooth formulation of plasticity with the Clausius-Duhem form of the second law of thermodynamics and with the energy balance (first law). Then, we note that the evolution of the entropy associated to the plastic deformation is nonsmooth.

\section{Nonsmooth mechanics and elastoplasticity}\label{nonsmooth_elastoplasticity}

We establish a link between variational multisymplectic formulation of continuum mechanics and elastoplasticity through internal slip of dislocation and internal friction due to lattice displacement (translation and rotation) which are dissipative nonsmooth dynamic phenomena.

\subsection{Moreau viewpoint and $1$D perfect plasticity} \label{Moreau}

Plastic bodies are characterized by the fact that their shape can be changed by the application of appropriately directed external forces, and that they retain their so-deformed shape upon removal of such forces. 

\medskip
%$\Delta_t\boldsymbol{\epsilon}_p(t,s)$ $\mathbf{v}_{\epsilon_p}(t,s)$

Let the \textit{internal stress $2$-tensor} $\boldsymbol{\sigma}(t,s)$ and $\left\llbracket \boldsymbol{\epsilon}_p \right\rrbracket_t$ a \textit{plastic strain-rate jump}\footnote{Recall that the \textit{jump} $\left\llbracket \cdot \right\rrbracket$ in \cite{Maugin1992} is defined by
\begin{equation}
\label{jump_n_form}
\left\llbracket V \right \rrbracket := V^+ - V^-
\end{equation}
where $V^-$ and $V^+$ are respectively referred to $t^-$ and $t^+$. In \cite{FeMaWe2003} the definition of the jump is extended to the spacetime.}, during time evolution,
 which occurs at position $s$ and time $t$.   
The \textit{yield criterion}\footnote{In this development we admit the existence of a single yield criterion in order to simplify the presentation, but generally there is a set of constraints.} $f:\boldsymbol{\sigma} \mapsto \mathbb{R}$, which confines the stresses $\boldsymbol{\sigma}$ to lie in the elastoplastic domain, is specified by the following inequality constraint 
\begin{equation}\label{sig_constraint}
f(\boldsymbol{\sigma})\leq 0.
\end{equation}
Thus the \textit{set of admissible stresses} is defined by 
\begin{equation}\label{elas_plas_dom}
\mathbb{E}_\sigma:  =  \{\boldsymbol{\sigma}  \, | \, f(\boldsymbol{\sigma})\leq 0 \} .
\end{equation}
The boundary of $\mathbb{E}_\sigma$ defined by $f(\boldsymbol{\sigma})=0$ is called the \textit{yield surface}. ``The points at which $\boldsymbol{\sigma}$ is inside the yield surface ($f(\boldsymbol{\sigma})<0$) constitute the elastic domain, while those where $\boldsymbol{\sigma}$ is on the yield surface form the plastic domain'' (\cite{Lubliner1990}).
The set $\mathbb{E}_\sigma$ is supposed to be convex. 
\medskip

Let $\mathbb{V}$ the \textit{set of plastic strain-rate jumps} $\left\llbracket \boldsymbol{\epsilon}_p \right\rrbracket_t$. The set of admissible stresses and the set of rate of change of plastic strain are placed in duality by a bilinear form $\left< \cdot, \cdot \right>$. 
The plasticity law is defined by stating the \textit{maximal dissipation principle}, i.e., the values of the stress $\boldsymbol{\sigma} \in \mathbb{E}_\sigma$ which correspond to some $\left\llbracket \boldsymbol{\epsilon}_p \right\rrbracket_t \in \mathbb{V}$ are the elements which minimize the numerical function $\boldsymbol{\sigma} \mapsto - \left<\left\llbracket \boldsymbol{\epsilon}_p \right\rrbracket_t, \boldsymbol{\sigma} \right>$.

\medskip

Next, we recall the convex analysis principles which allow to describe the plasticity.
Following \cite{Moreau1973, Moreau1976} the stress and the strain $(\boldsymbol{\sigma}, \left\llbracket \boldsymbol{\epsilon}_p \right\rrbracket_t)\in \mathbb{E}_\sigma \times \mathbb{V}$ which verify the plasticity law, i.e., the principle of maximum dissipation under the inequality constraint $f(\boldsymbol{\sigma})\leq 0$, can also be defined in an equivalent way by the \textit{variational inequality} with solution $\boldsymbol{\sigma}$ satisfying the condition
\begin{equation}
\left\{
\begin{aligned} \label{equiv_1}
\boldsymbol{\sigma} & \in \mathbb{E}_\sigma 
\\
\forall \boldsymbol{\sigma}' & \in \mathbb{E}_\sigma: \quad \left< \boldsymbol{\sigma} , \left\llbracket \boldsymbol{\epsilon}_p \right\rrbracket_t \right> \geq \left< \boldsymbol{\sigma}' ,\left\llbracket \boldsymbol{\epsilon}_p \right\rrbracket_t \right>,
\end{aligned}
\right.
\end{equation}
or in the following equivalent manners 
\begin{equation} \label{non_smooth_plast}
\begin{aligned}
\forall \boldsymbol{\sigma}'  \in \mathbb{R} &: \quad \left< \boldsymbol{\sigma}' - \boldsymbol{\sigma} ,\left\llbracket \boldsymbol{\epsilon}_p \right\rrbracket_t \right> + I_{\mathbb{E}_\sigma}(\boldsymbol{\sigma}) \leq   I_{\mathbb{E}_\sigma}(\boldsymbol{\sigma}'),
\\
\Leftrightarrow & \quad \left\llbracket \boldsymbol{\epsilon}_p \right\rrbracket_t \in \partial I_{\mathbb{E}_\sigma}(\boldsymbol{\sigma}),
\\
\Leftrightarrow & \quad \boldsymbol{\sigma} \in \partial I^*_{\mathbb{E}_\sigma}( \left\llbracket \boldsymbol{\epsilon}_p \right\rrbracket_t ),
\end{aligned}
\end{equation} 
where $I_{\mathbb{E}_\sigma}$ is the indicator function of $\mathbb{E}_\sigma$, i.e., $I_{\mathbb{E}_\sigma}(\boldsymbol{\sigma})=0$ if $\boldsymbol{\sigma} \in \mathbb{E}_\sigma$ and $I_{\mathbb{E}_\sigma}(\boldsymbol{\sigma})= + \infty$ if $\boldsymbol{\sigma} \notin \mathbb{E}_\sigma$. Its polar function $I^*_{\mathbb{E}_\sigma}$ is the support function\footnote{ See in \cite{RoWe1998} for a general development of the concept of \textit{support function}.} of $\mathbb{E}_\sigma$ relative to $\left< \boldsymbol{\sigma} , \left\llbracket \boldsymbol{\epsilon}_p \right\rrbracket_t \right>$, i.e.,
\begin{equation}\label{sup_function}
I^*_{\mathbb{E}_\sigma}(\left\llbracket \boldsymbol{\epsilon}_p \right\rrbracket_t) = \mathrm{sup}_{\boldsymbol{\sigma} \in T^2_0(\mathcal{S})} \left\{  \left< \boldsymbol{\sigma} , \left\llbracket \boldsymbol{\epsilon}_p \right\rrbracket_t \right> - I_{\mathbb{E}_\sigma}(\boldsymbol{\sigma}) \right\} = \mathrm{sup}_{\boldsymbol{\sigma} \in \mathbb{E}_\sigma} \left< \boldsymbol{\sigma} , \left\llbracket \boldsymbol{\epsilon}_p \right\rrbracket_t \right>.
\end{equation}
As a consequence \eqref{equiv_1} is also equivalent to 
\begin{equation}
\left\{
\begin{aligned} \label{equiv_2}
& \boldsymbol{\sigma}  \in \mathbb{E}_\sigma 
\\
 & \left< \boldsymbol{\sigma} , \left\llbracket \boldsymbol{\epsilon}_p \right\rrbracket_t \right>  = I^*_{\mathbb{E}_\sigma}(\left\llbracket \boldsymbol{\epsilon}_p \right\rrbracket_t ).
\end{aligned}
\right.
\end{equation}
That is, the values of $\boldsymbol{\sigma}$ associated with a given $\left\llbracket \boldsymbol{\epsilon}_p \right\rrbracket_t \in \mathbb{V}$, by the plasticity law, are the elements of $\mathbb{E}_\sigma$ for which the dissipated value is exactly equal to $I^*_{\mathbb{E}_\sigma}(\left\llbracket \boldsymbol{\epsilon}_p \right\rrbracket_t)$. In \cite{Moreau1973}, $I^*_{\mathbb{E}_\sigma}$ is denoted the \textit{dissipation function}. From now on we note the plastic dissipation by 
\begin{equation} \label{plast_dissip}
\mathcal{D}^p(\boldsymbol{\sigma}, \left\llbracket \boldsymbol{\epsilon}_p \right\rrbracket_t) := I^*_{\mathbb{E}_\sigma}(\left\llbracket \boldsymbol{\epsilon}_p \right\rrbracket_t).
\end{equation} 

 \begin{remark} {\rm
We recall that $\partial I_{\mathbb{E}_\sigma}(\boldsymbol{\sigma})=N_{\mathbb{E}_\sigma}(\boldsymbol{\sigma}) $ is the \textit{normal cone} to $\mathbb{E}_\sigma$ at $\boldsymbol{\sigma}$ and $\partial I^*_{\mathbb{E}_\sigma}(\left\llbracket \boldsymbol{\epsilon}_p \right\rrbracket_t) = N_{\mathbb{E}_\sigma}^*(\left\llbracket \boldsymbol{\epsilon}_p \right\rrbracket_t)$ is its dual; see \cite{RoWe1998}. Moreover we recall that $ \lambda \nabla_{\boldsymbol{\sigma}}f(\boldsymbol{\sigma}) \in N_{\mathbb{E}_\sigma}(\boldsymbol{\sigma})$, with the Lagrange multiplier $\lambda$ which satisfies $\lambda =0$ when $f(\boldsymbol{\sigma})<0$ and $\lambda\geq 0$ when $f(\boldsymbol{\sigma})=0$. Note that $\lambda$ obeys to the \textit{Kuhn-Tucker complementarity conditions}:
\[
\lambda \geq 0; \qquad f(\boldsymbol{\sigma})\leq 0; \qquad \lambda f(\boldsymbol{\sigma}) =0.  \qquad \square
\]}
\end{remark}

\begin{remark}Ê{\rm
An essential feature of the convex analysis principles is the impossibility to define $\left\llbracket \boldsymbol{\epsilon}_p \right\rrbracket_t$ as a single valued fonction of $\boldsymbol{\sigma}$, nor $\boldsymbol{\sigma}$ as a single valued function of $\left\llbracket \boldsymbol{\epsilon}_p \right\rrbracket_t$. Indeed for $\left\llbracket \boldsymbol{\epsilon}_p \right\rrbracket_t =0$ corresponds for $\boldsymbol{\sigma}$ all the values of $\mathrm{int}(\mathbb{E}_\sigma)$, and for $\boldsymbol{\sigma} \in \partial \mathbb{E}_\sigma$ corresponds for $\left\llbracket \boldsymbol{\epsilon}_p \right\rrbracket_t$ all the values of the normal cone $N_{\mathbb{E}_\sigma}(\boldsymbol{\sigma})$. Later on, we will tackle this issue. \quad $\square$}
\end{remark}

\begin{remark} {\rm
Concerning plasticity, it is important to note that the constraint \eqref{sig_constraint} is applied to the stress. There are no direct constraints on the configuration of the body. \quad $\square$ }
\end{remark}

\subsection{Vector measure and locally bounded variations} \label{blv_vm}

The vector measure plays an important role in the subsequent development, so we shall devote this section to recall the following results.

\paragraph{Vector measure.} Let a Banach space $X$, a locally compact domain $T$, and the vector space $\mathcal{H}(T)$ of real continuous functions $\psi:T \rightarrow \mathbb{R}$ with compact support. A \textit{vector measure} on $T$ with value in $X$, in the sense of \cite{Bourbaki1959}, is the linear application $m: \mathcal{H}(T) \rightarrow X$ such that for compact subdomain $K$ of $T$, the restriction of $m$ to $\mathcal{H}(T,K)$ is continuous for the topology of uniform convergence. If $\psi \in \mathcal{H}(T)$, the vector measure is noted $\int \psi \, dm$ instead of $m(\psi)$.

\medskip

Following \cite{Moreau1988b}, instead of a locally compact domain $T$, for simplification, we consider a real interval $I$ and we admit that $X$ is a Banach space, with metric denoted $\mathrm{d}$. Let $f:I \rightarrow X$ which is said to be of \textit{locally bounded variations} on $[a,b]$ iff $\mathrm{var}(f,[a,b]) < + \infty$; notation $f \in \mathrm{lbv}([a,b],X)$. Where the variation $\mathrm{var}(f,[a,b])$ of $f$ on $[a,b]$ is defined as follows
\begin{equation*}
\mathrm{var}(f,[a,b]) = \mathrm{sup} \sum_{i=1}^n \mathrm{d}\left(f(\tau^{i-1}), f(\tau^i)\right), \quad \text{with} \quad \tau^0=a, ..., \tau^n=b. 
\end{equation*}
In this definition the supremum is taken over all strictly increasing finite sequence $S: \, \tau^0<\tau^1< ... < \tau^n$ of points of $[a,b]$. 

From \cite{Moreau1988b} we recall the following results.
\begin{proposition} 
Let $f\in \mathrm{lbv}(I,X)$; for every $\psi \in \mathcal{H}(I)$ and every $\theta_S^i \in [\tau^{i-1},\tau^i]$, the mapping $S \mapsto  \sum_{i=1}^n \psi(\theta_S^i)(f(\tau^i) - f(\tau^{i-1}))$ converges to a limit independent of $\theta$, denoted $\int \psi df$. The convergence is uniform with regard to the choice of $\theta$. 
\end{proposition}

Note that the linear mapping $\mathcal{H}(I) \ni \psi \mapsto \int \psi df \in X$ constitues a vector measure on $I$ in the sense of Bourbaki. Where $df$ is called the \textit{differential measure (or Stieltjes measure)} of $f \in \mathrm{lbv}(I,X)$.

\paragraph{Radon-Nikodym property.} The Banach space $X$ has the \textit{Radon-Nikodym property} if, for every absolutely continuous $f:I \rightarrow X$, the differential measure $df$ admits a density $f'_t \in L^1(I,dt;X)$ relative to Lebesgue's measure $dt$; notation $df=f'_t dt$. Where $L^1(I,dt;X)$ is the notation for the set of $X$-valued functions which are $\mu$-integrable (in the sense of Bourbaki) over every compact subset of $I$.
In particular, the finite dimensional Banach space has the \textit{Radon-Nikodym property.}

An important result, when $f\in \mathrm{lbv}(I,X)$, deals with the $df$-measure of the singleton $\{a\}$. That is, for every $a\in I$ we have 
\begin{equation}\label{singleton_jump}
\int_{\{a\}} f'_t\,dt = f^+(a)-f^-(a).
\end{equation}
Thus, we deduce a relationship between \eqref{singleton_jump} and \eqref{jump_n_form} when the jump is locally bounded.

\begin{remark} {\rm
If $X$ has finite dimension, any $X$-valued measure $f$ is \textit{majorable}. That is, there exists a nonnegative real measure $\mu$ on $I$ such that, for every $\psi \in \mathcal{H}(I)$ one has $\|\int \psi \, df\| \leq \int | \psi |d \mu$.

Then it can be proved that, if $X$ is a finite dimensional Banach space, every $X$-valued vector measure $f$ possesses a density $f'_\mu \in L^\infty(I,\mu; X)$ relative to its modulus measure $\mu=|f|$. \qquad $\square$ }
\end{remark}

Another important result is the following
\begin{proposition} \label{derivative_density}
If $X$ possesses the Radom-Nikodym property and if $f\in \mathrm{lbv}(I,X)$, at Lebesgue almost every point $\tau$ of $I$, the function $f$ possesses a derivative function $\dot{f}(\tau)$ and after arbitrary extension to the whole of $I$, it constitutes a representative of the density $f'_t \in  L_{\mathrm{loc}}^1(I,dt;X)$. 
\end{proposition}

From this Proposition and \eqref{singleton_jump} we establish the link between the plastic strain-rate jump and a time derivative. Thus we get dimensionally consistent results as we will see. 

\begin{remark}{\rm
A natural generalization of the previous recapitulation dedicated to the vector measures consists in replacing the Lebesgue measure $dt$ by some prescribed nonnegative real measure $\mu$ on the interval $I$. \quad $\square$ }
\end{remark}

%\begin{remark} \label{scalar_measure} {\rm
\paragraph{Integral with respect to the vector measure.} Let the dual Banach space $X'$ of $X$ with duality pairing $\left<\cdot, \cdot \right>$. Given a vector measure $f$ on $I$ with value in $X$. For all $x' \in X'$, the linear mapping 
\begin{equation}\label{real_measure}
\mathcal{H}(I) \ni \phi \mapsto \left<x', f(\phi) \right> = \left<x', \int \phi \, df \right> =  \int \phi \, d(x'\circ f) \in \mathbb{R}
\end{equation}
is a real measure which depends linearly from $x'$. 

If $df=f'_\mu d\mu$, where $m'_\mu$ is the density with respect to the positive measure $\mu$ on $I$, the \textit{integral of $\phi$ with respect to the vector measure $f$} is defined by
\begin{equation}\label{integral_wr_m}
x'\circ f(\phi) = \int \phi \left< x', f'_\mu \right> d\mu,
\end{equation}
where $x'\circ f=\left< x', f'_\mu \right> \mu$ is a scalar measure for all $x' \in X'$, see \cite{Bourbaki1959}.

Later on, we will make some connections between the scalar mesure \eqref{integral_wr_m} and the dissipation \eqref{plast_dissip} previously defined.

\subsection{Variational formulation of nonsmooth mechanics} \label{FeMaWe}

\paragraph{Nonsmooth mechanics.} The fundamental theorem \ref{fundamental_th} which describe the variational multisymplectic formulation of nonsmooth continuum mechanics was presented in \cite{FeMaWe2003}. The general framework in which this theory was defined is field theory. The \textit{physical fields} $\varphi:X \rightarrow Y$ are the sections of the \textit{covariant configuration bundle} $\pi_{XY}:Y\rightarrow X$, where $X$ is the \textit{spacetime domain} with coordinates $\{x^0=t,\, x^1,...,x^n\}$ and $Y=X\times M$ is the \textit{configuration bundle} with \textit{ambient space} $M$ and coordinates $\{\varphi^1,..., \varphi^N\}$ on it. So the coordinates on $Y$ are written as $(x^\mu, \varphi^A)$ with $\mu=0,...,n$ and $A=1,...,N$.

The Lagrangian density is of the form
\begin{equation}\label{Lag_dens}
\mathcal{L} ( x^\mu, \varphi^A ,\dot{\varphi}^A , \varphi^A {}_{,i} )=  L (x^\mu, \varphi^A ,\dot{\varphi}^A , \varphi^A {}_{,i} )d^{n+1}x,
 \end{equation}
 where $\dot{\varphi}: = 
\partial \varphi/\partial t$ is the time derivative,$\varphi_{,i}: = \partial \varphi /\partial x^i$, $i=1,...,n$ are the partial space derivative, and $d^{n+1}x = dx^1\wedge ...\wedge dx^n\wedge dx^0$.
The associated action functional is defined to be
\begin{equation} \label{action_functional}
\mathfrak{S}^{ns} (\varphi):=\int_{U_X}
\mathcal{L} (x^\mu, \varphi^A , \dot{\varphi}^A , \varphi^A {}_{,i}).
\end{equation}
Stationarity of the action $ \mathfrak{S}^{ns} $ with respect to variations $ \delta \varphi $ yields the \textit{Euler-Lagrange field equations} or \textit{covariant Euler-Lagrange (CEL) equations} 
\begin{equation}\label{cEL}
\frac{\partial }{\partial t}   \frac{\partial L }{\partial \dot{\varphi}^A }+ \frac{\partial }{\partial x^i}   \frac{\partial L }{\partial \varphi_{,i}^A}- \frac{\partial L }{\partial\varphi^A }=0.
\end{equation}

\medskip

We introduce (see Figure below) a manifold $U$ with smooth closed boundary, a map $\phi: U\rightarrow Y$ taken to be smooth, the diffeomorphism $ \phi_X:U\rightarrow U_X\subset X$, and a submanifold $D \subset U$, called the \textit{singularity submanifold}, across which the Lagrangian $L$ may have singularities.
%\footnote{Contrary to what has been said in \cite{DeGBRa2016} this is the Lagrangian which may have singularities and not inevitably the field $\phi$. Indeed the singularities of $L$ can also be due to nonsmoothness of the first Piola-Kirchhoff stress tensor $\mathbf{P}(\varphi)$, or for other reasons.} 
Given $D_X := \phi_X(D)$, it is assumed that the singularity submanifold $D_X$ separates the interior of $U_X$ in two disjoint open subsets $U_X^+$ and $U_X^-$. The Lagrangian $L$ is assumed to be smooth only on $U_X\backslash D_X$. 
\begin{figure}[h!] {\footnotesize 
\begin{displaymath}
\begin{xy}
\xymatrix{   &    Y:= X \times M    \ar@<2pt>[d]^\pi  \\ 
U  \ar[ru]^{\phi}  \ar[r]^-{\phi_X}  &    U_X \subset X \ar@<2pt>[u]^{\varphi}  }
\end{xy}
\end{displaymath} }
%\caption{{\footnotesize Vector fields $V$ and $j^1V$ on $\pi :Y\rightarrow X$ and on $J^1Y\rightarrow X$.}}
\end{figure}

From the variational principle, as proved in \cite{FeMaWe2003}, we derive directly the equations of motion and the \textit{jump conditions}, staying on the Lagrangian side. In particular the jump conditions are due to the different orientations of $D_X$ when Stokes' theorem is involved in the integration by parts.

\begin{theorem} \label{fundamental_th}
Given a Lagrangian density $ \mathcal{L} ( x^\mu, \varphi^A ,\dot{\varphi}^A , \varphi^A {}_{,i} )$, which is smooth away from the discontinuity in $D_X$, there exists unique derivative of the action $\mathbf{d} \mathfrak{S}^{ns} (\varphi)$ such that for any $V =(V^\mu, V^A) \in T_\varphi \mathcal{C}$ compactly supported in $U$ and any open subset $U_X$ such that $\overline{U}_X \cap \partial X= \emptyset$,\footnote{Where $\overline{U}_X$ is the closure of $U_X$.}
\begin{align} 
\mathbf{d} \mathfrak{S}^{ns}( \phi ) \cdot (V) =&
\int_{U_X^+\cup U_X^-} \left( \frac{\partial L }{\partial\varphi^A } - \frac{\partial }{\partial t}   \frac{\partial L }{\partial \dot\varphi^A } - \frac{\partial }{\partial x^i}   \frac{\partial L }{\partial \varphi^A {}_{,i}} \right) \cdot V^A d^{n+1}x \label{nonsmooth_VP_local1}
\\
& + \int_{U_X^+\cup U_X^-}\left( \frac{\partial L}{\partial x^\nu }+ \frac{d}{d x^\mu} \left( \frac{\partial L}{\partial \varphi^A{}_{,\mu} } \varphi^A{}_{,\nu} \right) - \frac{d L}{d x^\nu}  \right)  V^\nu d^{n+1}x \label{nonsmooth_VP_local2}
\\
&+\int_{\partial U_X\backslash D_X}  \left(\frac{\partial L}{\partial \varphi^A {}_{, \mu } } \cdot V^A d^n x _\mu \right) \label{nonsmooth_VP_boundary1}
\\
&  + \int_{\partial U_X \backslash D_X}  \left(\left( L \delta ^\mu _\nu - \frac{\partial L}{\partial \varphi ^A {}_{, \mu } }\varphi ^A _{, \nu } \right)   V ^\nu   d^n x _\mu \right) \label{nonsmooth_VP_boundary2}
\\
&+ \int_{D_X }\left\llbracket \frac{\partial L}{\partial \varphi^A{}_{,\mu} } \cdot V^A  d^n x_\mu \right\rrbracket \label{nonsmooth_VP_local3}
\\
& +\int_{D_X }\left\llbracket \left( L\delta_\nu^\mu  - \frac{\partial L}{\partial  \varphi^A{}_{,\mu} } \varphi^A{}_{,\nu} \right)  V^\nu  d^n x_\mu \right\rrbracket. \label{nonsmooth_VP_local4}
\end{align}

\end{theorem}

Where \eqref{nonsmooth_VP_local1} gives the CEL equations. The time component of \eqref{nonsmooth_VP_local2} is the \textit{energy-evolution equation}, while its full expression formulates the \textit{balance of configurational forces}.
The last two expressions \eqref{nonsmooth_VP_local3} and \eqref{nonsmooth_VP_local4} are respectively the \textit{vertical jump conditions} involving momenta and the \textit{horizontal jump conditions}, i.e., energy jump conditions, which are the consequence of local nonsmoothness when $x\in D_X$.

\begin{remark}\label{Radon_Niko} {\rm
Let an application $U_X \ni x\mapsto V(x) \in \mathbb{R}^N$ with locally bounded variations on $U_X$, see in \S \ref{blv_vm} for the definition. Where this function may have discontinuities but, at every point $x$, the left-limit and the right-limit exist. Then one can associate an $\mathbb{R}^N$-valued measure. If one has $V^-(x) \neq V^+(x)$ for $x\in D_X$, where $V^-$ and $V^+$ are respectively referred to $U_X^-$ and $U_X^+$, as explained in \cite{Moreau1988a} the $\mathbb{R}^N$-valued measure possesses one atom $x\in D_X$ with jump value $\left\llbracket V \right \rrbracket (x) =V^+(x)- V^-(x)$.

 Also the jump can be interpreted as the value of a \textit{vector measure}, respectively \textit{covector measure}, on the set $ D_X$ of the atoms. 

Let the case, where we have a form $\boldsymbol{\alpha}$. In view of the Radon-Nikodym property, every covector measure $\left\llbracket \boldsymbol{\alpha} \right \rrbracket (x)$ on the locally compact subset $D_X$ of $U_X\subset \mathbb{R}^{n+1}$ may be represented as follows: there exists a positive scalar measure $dt$ on $D_X$ and a \textit{density function} denoted $\left\llbracket \boldsymbol{\alpha} \right \rrbracket_{t}'$ relatively to the Lebesgue measure, such that one writes
\begin{equation} \label{jump_vector_measure}
j^*\left\llbracket \boldsymbol{\alpha} \right \rrbracket (x)= \left\llbracket \boldsymbol{\alpha} \right \rrbracket_{t}'(x) dv_{D_X},
\end{equation}
where $j: D_X \hookrightarrow U_X$ is the inclusion and $dv_{D_X}$ is the measure associated to the volume form on $D_X$. 

However even though we know from the Radon-Nikodym theorem that there exists a density function, the theorem does not indicate how to calculate this density function. This question will be solved through convex analysis.
\quad $\square$
}
\end{remark}

\paragraph{Nonsmooth continuum mechanics and inequality constraints.} The variational inequalities and the problems of constrained minimization were widely studied by Moreau and Rockafellar, see e.g., \cite{RoWe1998}. From this viewpoint there have been developments that bear on the multisymplectic formulations of nonsmooth continuum mechanics when the configuration is subjected to inequality constraints, see e.g., \cite{FeMaOrWe2003} and \cite{DeGBRa2016}.

\medskip

For example, for frictionless contacts, the force of constraint is normal to the concerned bodies. Following \cite{Moreau1988a}, the $dt$-measurable vector field $\left\llbracket \boldsymbol{\alpha} \right \rrbracket_{t}'$ satisfies   
\begin{equation}\label{vert_jump_vect_meas}
 \left\llbracket \boldsymbol{\alpha} \right \rrbracket_{t}'(x) \in \partial I_C(\varphi(x)),
 \end{equation}
where $\partial I_C(\varphi(x))$ is the normal cone\footnote{We recall that the vectors in the normal cone $\partial I_C(\varphi(x))$ are of the form $ \boldsymbol{\lambda} (x) \cdot\nabla \psi( \varphi (x))$ where $\boldsymbol{\lambda}$ are the Lagrange multipliers and $\psi$ are the inequalities constraints $\psi(\varphi)\leq 0$.} to the admissible contact domain $C\subset M$ in $\varphi(x)$ and the density function denoted
\[
  \left\llbracket  
\frac{\partial L}{\partial \varphi^A {}_{, \nu } }(x) N_\nu (x)
\right\rrbracket_t'
\] 
is defined to be the vector field associated to the vertical jump \eqref{nonsmooth_VP_local3}, where  $N_\nu$ is the normal vector to $D_X$.

From \eqref{vert_jump_vect_meas} and the properties of a normal cones we can deduce that there exists a Lagrange multiplier $\boldsymbol{\lambda}$ such that we get for every $V^A$  
%, and the duality way of defining a vector measure through the pairing (see \cite{Bourbaki1959}),
\begin{equation}\label{contact_vert_jump}
\int_{D_X }\left\llbracket  
\frac{\partial L}{\partial \varphi^A {}_{, \nu } }(x) N_\nu (x)
\right\rrbracket_t' \!\! \cdot V ^A(\varphi(x))  dv_{D_X} = \int_{D_X}  \lambda_i (x)\psi^i_{,A}( \varphi (x))\cdot V ^A(\varphi(x)) dv_{D_X} . 
\end{equation}
  If we consider the time component relative to the Lebesgue's measure $dt$, due to Proposition \ref{derivative_density}, we can deduce that there exists a time derivative which constitutes a representative of the momenta jump. As a confirmation, the equation \eqref{contact_vert_jump} is dimensionally consistent, i.e., on the left we have the time derivative of a momenta and a force on the right.  
 Concerning the left-hand side of \eqref{contact_vert_jump}, one has 
\begin{equation*} %\label{contact_vert_jump_2}
\int_{D_X }\left\llbracket  
\frac{\partial L}{\partial \varphi^A {}_{, \nu } }(x) N_\nu (x)
\right\rrbracket_t' \!\! \cdot V ^A(\varphi(x))  dv_{D_X}  = \int_{D_X }\left\llbracket  
\frac{\partial L}{\partial \varphi^A {}_{, \nu } }(x)N_\nu(x) \cdot V^A(\varphi(x))   
\right\rrbracket_t'  dv_{D_X},
\end{equation*}
where we used the continuity of the vector field $V$. Also we get the statement established in \cite{DeGBRa2016}.

\subsection{Perfect elastoplasticity and nonsmooth mechanics} \label{non_smooth_elastoplastic}

We refer to \cite{Moreau1973,Moreau1976,Moreau1982,Moreau1988a,Moreau1988b} in addition to \cite{FeMaWe2003} for this development. 
%Where the main point is as follows. 

In order to associate the variational multisymplectic formulation and convex analysis, we will take into account the variational inequalities \eqref{equiv_1} where constraints \eqref{sig_constraint} are included. Additionally, we will use the properties of the variational inequalities which can be expressed in different equivalent ways, in particular through the plastic dissipation.

\subsubsection{Lagrangian density and plastic dissipation} 

\paragraph{Lagrangian density.} The \textit{material frame indifference} states that if we view the configuration from a rotated point of view, then the stress transforms by the same rotation. Also, if we want that the stored energy function $W(x, \varphi , \mathbf{F}_e)$ satisfies this property it must depends on the elastic gradient deformation $\mathbf{F}_e$ through right Cauchy-Green deformation tensor $\mathbf{C}_e=\mathbf{F}_e^T \mathbf{F}_e$ (see \cite{MaHu1994}).

We admit that the deformation gradient takes the form of a local multiplicative decomposition into elastic and plastic matrices as $\mathbf{F}=\mathbf{F}_e\mathbf{F}_p$. 

The Lagrangian density is defined as follows
\begin{equation} \label{Lagrangian_eucl}
 \mathcal{L} ( x, \varphi ,\dot{\varphi} , \mathbf{F}_e) = \frac{1}{2} \rho(x) \left<\dot{\varphi}, \dot{\varphi}\right> dv(x) - \rho(x) W(x, \varphi , \mathbf{C}_e( \mathbf{F}_e))dv(x), 
\end{equation}
with the mass density $\rho$ and the elastic component $\mathbf{F}_e =\mathbf{F}\mathbf{F}_p^{-1}$ of the deformation gradient. The components of $\mathbf{F}$ are $F^A{}_{i}=\varphi^A {}_{,i}$. Note that, for simplicity, we consider the Euclidean case.

\paragraph{Plastic dissipation through viscous regularization.} The problem to resolve is to calculate the plastic strain-rate jump and the plastic dissipation, while the material frame indifference is required. The difficulty is precisely that the plasticity is a nonsmooth phenomena, see \S \ref{Moreau}, where the time rate of change of the plastic components $\mathbf{F}_p$ evolve by jumps.
%, denoted $\left\llbracket \mathbf{F}_p \right\rrbracket_t $.

To solve this difficulty we recall from \cite{Moreau1973} that we can take into account at the same time several
resistance laws, like viscosity and plasticity. The introduction of viscosity produces a regularization effect of the plastic strain called \textit{Moreau-Yosida regularization}. Conversely, the plasticity can be seen as the limit of the viscoplasticity when viscosity disappears. Indeed, adding a tiny viscosity effect to a plasticity law produce a penalty function, where the size of the penalty coefficient is inversely proportional to the value of the viscosity coefficient $\eta$. But, when $\eta \rightarrow 0$, the penalty function becomes the indicator function. So we get again a plastic law and a strain-rate which evolves by jumps. 

The viscoelastoplasticity is characterized by $ - \mathbf{S} : \mathbf{D}_p + \eta \, \mathbf{q}(d(\mathbf{S}))$, which correspond to the viscous regularization of the plastic dissipation $\mathbf{S} : \mathbf{D}_p$ through the viscous dissipation $\eta \, \mathbf{q}(d(\mathbf{S}))$. Where $\mathbf{S}$ is the symmetric second Piola-Kirchhoff stress tensor, $\mathbf{D}_p$ is the plastic strain-rate, $\mathbf{q}$ is a quadratic form, and $d(\mathbf{S})$ is the distance between $\mathbf{S}$ and the convex admissible elastoplastic domain $\mathbb{E}_{\mathbf{S}}$, characterized by the yield conditions $ f(\mathbf{S})\leq 0$, i.e., 
  \begin{equation}\label{admis_Piola_stress}
\mathbb{E}_{\mathbf{S}}:=  \{\mathbf{S}  \, | \, f(\mathbf{S})\leq 0 \}.
\end{equation}
Thus the Moreau-Yosida regularization occurs when $ f(\mathbf{S})\geq 0$. However, as before, we have an elastic phenomena when $f(\mathbf{S})<0$. 

Due to the regularization effect, induced by the introduction of viscosity, we can calculate $\mathbf{D}_p$  %, see \cite{Simo1998}\,(36.11) 
through the derivative $\dot{ \mathbf{F}}_p$ (see in \cite{Simo1998} for details)
 \begin{equation}\label{Second_Pio_Kirch}
 \mathbf{D}_p:= \left( \mathbf{C}_e  \dot{\mathbf{F}}_p  \mathbf{F}_p^{-1} \right)_{\mathrm{sym}}
 \end{equation}
 where $\mathbf{C}_e$ and $\dot{\mathbf{F}}_p  \mathbf{F}_p^{-1}$ are unaffected by rigid motions superposed on the current placement. 

In order to keep only the plastic law, we remove the viscosity ($\eta=0$). Hence the elastoplastic strain-rate jump, denoted $ \left\llbracket \mathbf{D}_p\right\rrbracket_t$, and the elastoplastic dissipation $\mathcal{D}^p$ are defined as the following limits
 \begin{equation}\label{plastic_strain_rate}
 \begin{aligned}
 & \mathrm{lim} \,{}_{\eta \rightarrow 0} \left( \mathbf{C}_e  \dot{\mathbf{F}}_p  \mathbf{F}_p^{-1} \right)_{\mathrm{sym}}  =:  \left\llbracket \mathbf{D}_p\right\rrbracket_t  ,  
  \\
& \mathrm{lim} \,{}_{\eta \rightarrow 0} \;\mathbf{S} : \mathbf{D}_p =\; \mathbf{S} : \left\llbracket \mathbf{D}_p\right\rrbracket_t =  \mathcal{D}^p  \; ,
  \end{aligned}
 \end{equation}
 where the values of $\mathbf{C}_e$ is given at time $t^-$. From \eqref{non_smooth_plast} we recall that 
 \begin{equation}\label{non_smooth_plast_2}
 \left\llbracket \mathbf{D}_p\right\rrbracket_t \in \partial \mathrm{I}_{\mathbb{E}_{\mathbf{S}}}(\mathbf{S}).
 \end{equation}
 
 Note that the rate of plastic deformation $\left\llbracket \mathbf{D}_p\right\rrbracket_t$ is different from zero if and only if $f(\mathbf{S})= 0$, i.e., the plastic components $\mathbf{F}_p$ are preserved when $f(\mathbf{S})\neq 0$. 
 %Indeed the normal cone $\partial \mathrm{I}_{\mathbb{E}_{\mathbf{S}}}(\mathbf{S})$ to $\mathbb{E}_{\mathbf{S}}$ in $\mathbf{S}$ is equal to $\{0\}$ when $\mathbf{S} \in \mathrm{int}(\mathbb{E}_{\mathbf{S}})$.
In addition, we recall that the normal cone $\partial \mathrm{I}_{\mathbb{E}_{\mathbf{S}}}(\mathbf{S})$ in $\mathbf{S}$ can be written in the form $\boldsymbol{\lambda} \partial_{\mathbf{S}} f(\mathbf{S})$, where $\boldsymbol{\lambda}$ are Lagrange multipliers.
 
\begin{remark} {\rm
It is important to note that the nonsmoothness of the plastic gradient of deformation $\mathbf{F}_p$, which evolves by jumps, induces nonsmoothness of the density Lagrangian \eqref{Lagrangian_eucl}. So, in the following, we can take into account of the Theorem \ref{fundamental_th}. \quad $\square$ }
\end{remark}

\begin{remark}{\rm 
The introduction of viscosity in elastoplasticity and the Moreau-Yosida regularization are studied more in details in \cite{Demoures2018b}.\quad $\square$
 }
 \end{remark}

\subsubsection{Variational formulation of perfect elastoplasticity}

We know that outside of the nonsmooth plastic phenomena the elastic deformation are described by the CEL equations \eqref{cEL}. Also we will investigate in the multisymplectic framework the vertical and horizontal jumps. 

\medskip

\paragraph{Singularity submanifold.} The submanifold $D_X\subset U_X$ is the singularity submanifold which matches with the plastic domain $\partial \mathbb{E}_{\mathbf{S}}$ defined by $f(\mathbf{S}(\varphi))=0$. That is $D_X$ is the space-time domain locally compact such that the time evolution of $\mathbf{S}(\varphi)$ is nonsmooth and the dissipation function $\mathcal{D}^p$, defined in \eqref{plastic_strain_rate}, verifies the maximal dissipation principle under inequality constraint $f(\mathbf{S})\leq 0$.  We can deduce that exists a Lagrangian multipliar $\boldsymbol{\lambda}$ such that $\mathcal{D}^p$ ascertains the maximal dissipation when $\left\llbracket \mathbf{D}_p\right\rrbracket_t = \boldsymbol{\lambda} \partial_{\mathbf{S}} f(\mathbf{S})$. 

%From now on we denote by $\overline{\mathbf{S}}:= \mathbf{S}(\varphi(x))$ when $ x \in D_X$. We recall that it verify $\overline{\mathbf{S}} \in \partial I^*_{\mathbb{E}_{\mathbf{S}}}( \left\llbracket \mathbf{D}_p\right\rrbracket_t )$ (see in \S \ref{Moreau}).

%\textbf{1)} 
\paragraph{Vertical variations ($V^A \neq 0$).} The plasticity is due to slips and frictions inside the body, but without impenetrability constraints on the fields $\varphi$, contrary to what is happening in contact mechanics, see, e.g., \eqref{contact_vert_jump}. So we get from \eqref{contact_vert_jump} the following jumps conditions without reaction forces (contact impulsion)
\begin{equation}\label{vertical_jump}
 \int_{D_X }\left\llbracket \frac{\partial L}{\partial \varphi^A{}_{,\nu} }(x) N_\nu(x) V^A \right\rrbracket_t' \!    dv_{D_X}  = 0. 
\end{equation}

\paragraph{Horizontal variations ($V^\mu \neq 0$).}   
 
Given the fact that the set of plastic deformation $2$-tensor is a finite dimensional Banach space with the Radon-Nikodym property, the rate of plastic deformation $\left\llbracket \mathbf{D}_p\right\rrbracket_t$ is the value of a vector measure with bounded variations which possesses a density denoted $\left\llbracket \mathbf{D}_p\right\rrbracket_t'$ relative to a measure $dt$, see in \S \ref{blv_vm}. Moreover the dissipation $\mathbf{S} : \left\llbracket \mathbf{D}_p\right\rrbracket_t$ is a pairing between the stress and the rate of plastic deformation. Such that we get a scalar measure $\left(\mathbf{S} : \left\llbracket \mathbf{D}_p\right\rrbracket_t' \right) dt$ for all $\mathbf{S}$ (see \eqref{integral_wr_m}).
In addition, from Proposition \ref{derivative_density} we deduce that the units are respected, i.e., the elastoplastic dissipation $\mathcal{D}^p$, as defined in \eqref{plastic_strain_rate}, is a power. 

Concerning the time component of the horizontal jump \eqref{nonsmooth_VP_local4}, i.e., $V^0\neq 0$, $V^j=0$ for $j=1, ..., n$. This is the value of a real measure with bounded variations which possesses a density denoted  $\left\llbracket  L  - \frac{\partial L}{\partial  \varphi^A{}_{,\mu} } \dot{\varphi}^A \right\rrbracket_t'$ relative to a measure $dt$, which can be represented by a time derivative.
As a consequence the global energy jump condition is given with the correct units, by
\begin{equation}\label{horiz_jump}
\int_{D_X} \left\llbracket  L  - \frac{\partial L}{\partial  \varphi^A{}_{,\mu} } \dot{\varphi}^A  \right\rrbracket_t' V^0 dv_{D_X} - \int_{D_X}  \left(\mathbf{S} : \left\llbracket \mathbf{D}_p\right\rrbracket_t' \right) V^0 dv_{D_X} \ni 0,
\end{equation}
which characterizes the intersection between the horizontal jump condition due to nonsmoothness and the maximum dissipation principle.

\medskip

The previous results leads to the following theorem which describes the elastoplastic behavior through the variational multisymplectic formulation.
\begin{theorem} \label{plasticity_vert_hor}
Consider a Lagrangian density $\mathcal{L} ( x^\mu, \varphi^A ,\dot{\varphi}^A  , F^A{}_{i} )$ which is smooth away from the discontinuity $D_X$, where deformation gradient $\mathbf{F}=\mathbf{F}_e\mathbf{F}_p$ is seen as the composition of elastic and plastic deformations. We assume the same hypotheses as above on 
$\pi _{XY}: Y \rightarrow X$. 
Then $ \phi=( \phi _{X}, \varphi )$ is a critical point of $\mathfrak{S}^{ns}$ 
relative to the constraint \eqref{admis_Piola_stress} on the second Piola-Kirchhoff stress tensor $\mathbf{S}$ if and only if 
\begin{itemize}
\item Away from the singularity, the field $\phi$ satisfies the covariant Euler-Lagrange equations on $U_X \backslash D_X$, 
\begin{equation}\label{CELeq}
\frac{\partial L }{\partial\varphi^A } - \frac{\partial }{\partial t}   \frac{\partial L }{\partial \dot\varphi^A } - \frac{\partial}{\partial x^i} \frac{\partial L }{\partial \varphi^A {}_{,i}} =0,
\end{equation}  % \frac{\partial L }{\partial \varphi^A {}_{,i}}
together with the balance of energy on $U_X \backslash D_X$. 

\item  When $x \in D_X$ the field $ \phi$ verify the following conditions:
\\ 
{\rm \textbf{(a)}} the vertical jump condition: 
\begin{equation}\label{vert_jump_cond}
\int_{D_X }\left\llbracket \frac{\partial L}{\partial \varphi^A{}_{,\nu} }(x) N_\nu(x) V^A \right\rrbracket_t' \!    dv_{D_X}  = 0. 
\end{equation}  
\\ 
{\rm \textbf{(b)}} the global energy jump condition (time component): for all vector 
fields $V^0$ we have
\begin{equation}\label{hor_jump_cond}
\int_{D_X} \left\llbracket L  - \frac{\partial L}{\partial  \varphi^A{}_{,\mu} } \dot{\varphi}^A  \right\rrbracket_t' V^0 dv_{D_X} = \int_{D_X}  \left(\mathbf{S} : \left\llbracket \mathbf{D}_p\right\rrbracket_t' \right) V^0 dv_{D_X}, 
\end{equation}
where $\left\llbracket \mathbf{D}_p\right\rrbracket_t = \boldsymbol{\lambda} \partial_{\mathbf{S}} f(\mathbf{S}) \in \partial I_{\mathbb{E}_P}( \mathbf{S}) $ with the Lagrange multiplier $\boldsymbol{\lambda}$. 
\item 
On the boundary $\partial U_X \setminus D_X$, the field $ \phi $ 
verifies the following conditions:\\
{\rm \textbf{(c)}} we have
\begin{equation}
\label{first_partial_bc}
\frac{\partial L}{\partial \varphi ^A {}_{, \nu }}  V ^A  N_\nu = 0,
\end{equation}
{\rm \textbf{(d)}} for all $\nu=0,...,n$ we have
\begin{equation}
\label{second_partial_bc}
\left( L \delta ^\mu _\nu - 
\frac{\partial L}{\partial \varphi ^A {}_{, \mu } }\varphi ^A{} _{, \nu } \right) N_\mu = 0.
\end{equation} 
\end{itemize}
\end{theorem}

\begin{remark}{\rm
``The fact that the constraints involve only spatial and not time derivatives means that imposing the constraints is equivalent to restricting the infinite-dimensional configuration manifold used to formulate the theory as a traditional Hamiltonian or Lagrangian field theory. In this case, the constraint is simply a holonomic or configuration constraint and it is known that restricting HamiltonÕs principle to the constraint submanifold gives the correct equations for the system. '' \cite{MaPeShWe2001} \quad $\square$ }
\end{remark}

\paragraph{Noether theorem.} Consider a one-parameter family $\phi^\epsilon$ of deformation mappings that are a symmetry of the mechanical system. That is the Lagrangian $\mathcal{L}$  is equivariant with respect to the symmetry group action $G$. This implies the preservation of the action functional $\mathfrak{S}^{ns}(\phi) =\int_{U_X} \mathcal{L} ( x, \varphi^A ,\dot{\varphi}^A , F^A{}_{i} ) $, i.e., $\mathfrak{S}^{ns}(\phi) = \mathfrak{S}^{ns}(\phi^\epsilon)$ where $\phi^0=\phi$.
When the plastic dissipation occurs we can deduce the integral form of Noether's theorem, from \cite{DeGBRa2016} in \S 2.4. 

Due to vertical jump conditions in \eqref{vert_jump_cond} we get:
For all open subsets $U' \subset U$ with piecewise smooth boundary and for all $\xi$ in the Lie algebra $\mathfrak{g}$ of the Lie group $G$, we have
\begin{equation}\label{Plastic_Noether} 
\int_{\phi _X(D')} \left\llbracket \frac{\partial L}{\partial \varphi^A{}_{,\mu} }  \xi_Y^A \, d^n x_\mu \right\rrbracket  + \int_{\partial \phi _X(U')\setminus \phi_X(D')} 
\left(\frac{\partial L}{\partial \varphi^A {}_{, \mu } }  \xi_Y^A \, d^n x _\mu \right) =0,
\end{equation} 
where $ \xi_Y $ is the infinitesimal generator.

\begin{remark}{\rm
Generally nonsmoothness is associated with boundary contact, and friction, and/or with interior plasticity. In such cases the question of conservation of symmetries must be studied under combination of different perspectives. \quad $\square$ }
\end{remark}

\paragraph{Multisymplectic form formula.} Concerning the Cartan form (or multisymplectic form) and the multisymplectic form formula this will be studied in its own for various types of nonsmooth problems in a paper to come.

\section{Rheological model for nonsmooth elastoplasticity} \label{rheo_nonsmooth_elasto}

The rheological model we consider was described in \cite{SiHu1998}.

\subsection{Elastoplastic model} \label{elas_plast_mod}

\begin{figure}[H] \centering 
 \includegraphics[width=2.25 in]{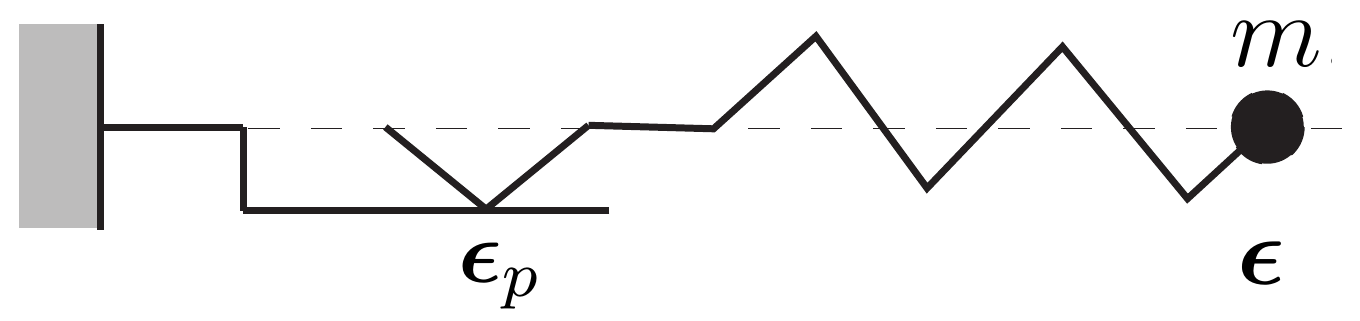} \vspace{-3pt}
 \caption{\footnotesize Elastoplastic bar.  } \label{elast_bar_0} 
 \end{figure}

The $1$D rheological model (see Fig \ref{elast_bar_0}) is composed of an elastic spring of length $\ell_0$ at rest, with Young modulus $E$, and elastic strain $\boldsymbol{\epsilon}_e$. 
At one end of the spring we fix a mass $m$, and at the opposite we fix a frictional pad which induces the frictional strain $\boldsymbol{\epsilon}_p$ referred to as the plastic strain. 
  We admit that the total strain denoted $\boldsymbol{\epsilon}$ is defined as the sum of the elastic strain $\boldsymbol{\epsilon}_e$ and the plastic strain $\boldsymbol{\epsilon}_p$, i.e., 
\begin{equation} \label{total_displ}
\boldsymbol{\epsilon}:= \boldsymbol{\epsilon}_e+ \boldsymbol{\epsilon}_p.
\end{equation}

The elastic strain is measured as the change in length $\Delta \ell$ divided by the original length $\ell_0$. 
Concerning our model we choose $\ell_0=1$. This hypothesis allows to identify the elongation of the spring with the elastic strain $\boldsymbol{\epsilon}_e$. At initial time $t^0$ we admit that the plastic strain $\boldsymbol{\epsilon}_p$ is $0$.
The decomposition \eqref{total_displ} is valid when plastic strain and elastic strain are small, see \cite{Lubliner1990} p.486 and \cite{Maugin2011} p.46. Which fits well with the rheological models.

We deduce from our model that the stress due to the elastic strain and the stored energy associated are respectively given by 
\begin{equation} \label{stress}
\boldsymbol{\sigma}= E \boldsymbol{\epsilon}_e = E(\boldsymbol{\epsilon}- \boldsymbol{\epsilon}_p) \quad \text{and} \quad W(\boldsymbol{\epsilon}) = \frac{1}{2} E |\boldsymbol{\epsilon}- \boldsymbol{\epsilon}_p|^2.
\end{equation}

 \subsection{Perfect plasticity case} \label{perfect_plastic}

The frictional device is characterized by the yield criterion $f$ which constrains the admissible stress $\boldsymbol{\sigma}$ to lie in the admissible set \eqref{elas_plas_dom}.

The Lagrangian $L(\boldsymbol{\epsilon}, \dot{\boldsymbol{\epsilon}})$ associated to perfect plasticity is defined by
\begin{equation}\label{Lagrangian_0}
\begin{aligned}
L(\boldsymbol{\epsilon}, \dot{\boldsymbol{\epsilon}})  & = \frac{1}{2} m |\dot{\boldsymbol{\epsilon}}|^2 - \frac{1}{2} E |\boldsymbol{\epsilon}_e|^2
\\
& = \frac{1}{2} m |\dot{\boldsymbol{\epsilon}}|^2 - \frac{1}{2} E |\boldsymbol{\epsilon}- \boldsymbol{\epsilon}_p|^2.
\end{aligned}
\end{equation}
Due to the plasticity laws described in \S \ref{Moreau}, the rate of change of total strain is given by $\dot{\boldsymbol{\epsilon}}=\dot{\boldsymbol{\epsilon}}_e + \left\llbracket \boldsymbol{\epsilon}_p \right\rrbracket_t$, where $\left\llbracket \boldsymbol{\epsilon}_p \right\rrbracket_t =0$ when $x\notin D_X$, i.e., $\boldsymbol{\epsilon}_p = Cte$ apart from plastic domain.

From Theorem \ref{plasticity_vert_hor}, we can describe the different situations that one meets: 

\paragraph{1. Elastic regime.} The elastic regime prevails as long as $\left\llbracket \boldsymbol{\epsilon}_p \right\rrbracket_t=0$. It is equivalent to saying that $\boldsymbol{\epsilon}_p=$ constant. The unconstrained CEL equations \eqref{EL} describe the motion of the system. Whenever the CEL equations are satisfied, 
the time energy-evolution equation is equal to zero. 
On the time interval $[0,T]$ the action map to be
\begin{equation}\label{action_0}
\mathfrak{S}^{ns}(\boldsymbol{\epsilon}) = \int_{0}^{T} L(\boldsymbol{\epsilon}, \dot{\boldsymbol{\epsilon}}) dt.
\end{equation}
Computing the variation of the action map $\mathbf{d}\mathfrak{S}^{ns}(\boldsymbol{\epsilon})\cdot \delta \boldsymbol{\epsilon}$ we get, from the Hamilton principle, the \textit{Euler-Lagrange equations}
\begin{equation}\label{EL}
m\ddot{\boldsymbol{\epsilon}} + E(\boldsymbol{\epsilon}- \boldsymbol{\epsilon}_p) =0.
\end{equation}
In addition from the horizontal variations $\delta t$ we get the conservation of energy.

\paragraph{2. Plastic regime.} The plastic regime occurs at time $\overline{t}$ when $f(\overline{\boldsymbol{\sigma}})=0$ holds
\begin{enumerate}
\item the laws of plasticity \eqref{non_smooth_plast} are verified, i.e., $\left\llbracket \boldsymbol{\epsilon}_p \right\rrbracket_t \in N_{\mathbb{E}_\sigma}(\overline{\boldsymbol{\sigma}})$,

\item From Theorem \ref{plasticity_vert_hor} it follows that
\begin{enumerate}
\item the vertical jump condition induces momenta conservation $\partial K(\dot{\boldsymbol{\epsilon}})/\partial \dot{\boldsymbol{\epsilon}}$ due to the absence of constraints on $\boldsymbol{\epsilon}$; see \eqref{vert_jump_pp}.
\item the horizontal energy jump $\left\llbracket -E_{\mathrm{tot}} \right\rrbracket$ is exactly equal to the plastic dissipation $\left< \overline{\boldsymbol{\sigma}}, \left\llbracket \boldsymbol{\epsilon}_p \right\rrbracket_t \right>$; see \eqref{rate_energ_0}.
\end{enumerate}
\end{enumerate}
From the vertical jump condition \eqref{vert_jump_cond} we get  
\begin{equation}\label{vert_jump_pp}
\left\llbracket m\dot{\boldsymbol{\epsilon}} \right\rrbracket  =0. 
\end{equation}
Due to plastic laws the rate of change of plastic strain satisfies \eqref{non_smooth_plast}, i.e.,
\begin{equation}\label{sigma_constraint}
\left\llbracket \boldsymbol{\epsilon}_p \right\rrbracket_t \in N_{\mathbb{E}_\sigma}(\overline{\boldsymbol{\sigma}}).
\end{equation}

At time $\overline{t}$ the horizontal jump conditions \eqref{hor_jump_cond} give the rate of change of the total energy, i.e., there exists $\lambda$ such that  
\begin{equation}\label{rate_energ_0}
\begin{aligned}
 \left\llbracket - E_{\rm tot}  \right\rrbracket &  =  \left< \overline{\boldsymbol{\sigma}} , \lambda \partial_{\boldsymbol{\sigma}} f(\overline{\boldsymbol{\sigma}}) \right> =  \mathcal{D}^p(\overline{\boldsymbol{\sigma}}, \lambda \partial_{\boldsymbol{\sigma}} f(\overline{\boldsymbol{\sigma}})),
\end{aligned}
\end{equation}
where
\begin{equation} \label{tot_energy_0}
 \textcolor{red}{\left\llbracket \boldsymbol{\epsilon}_p \right\rrbracket_{\bar{t}} }= \lambda \partial_{\boldsymbol{\sigma}} f(\overline{\boldsymbol{\sigma}}) \quad \text{and}\quad E_{\rm tot} (\boldsymbol{\epsilon}, \dot{\boldsymbol{\epsilon}}) =  \frac{1}{2} m |\dot{\boldsymbol{\epsilon}}|^2 + \frac{1}{2} E |\boldsymbol{\epsilon}- \boldsymbol{\epsilon}_p|^2.
\end{equation}

\begin{example} \label{exemple_1}
 Let the \textit{Tresca criterion} defined in \cite{Tresca1872} as follows
\begin{equation}\label{Tresca_law}
f(\boldsymbol{\sigma}) := |\boldsymbol{\sigma}| - \boldsymbol{\sigma}_Y \leqslant 0, \quad \text{with} \quad \boldsymbol{\sigma}_Y > 0.  
\end{equation}
We deduce the expression of the normal cone to $\mathbb{E}_\sigma$ at $\boldsymbol{\sigma}$
\begin{equation}
N_{\mathbb{E}_\sigma}(\boldsymbol{\sigma})= \lambda \, \partial f(\boldsymbol{\sigma} ) = \lambda \, \mathrm{sgn}(\boldsymbol{\sigma} ), \quad \text{for all} \quad \lambda \geqslant 0.  \nonumber %\label{norm_cone_0}
\end{equation}
The plastic phenomenon occurs when $f(\boldsymbol{\sigma})=0$, i.e., when $|\overline{\boldsymbol{\sigma}} | = \boldsymbol{\sigma}_Y$. We get
\begin{equation} \label{cond_fric_0}
\left\llbracket \boldsymbol{\epsilon}_p \right\rrbracket_{\bar{t}} = \lambda \, \mathrm{sgn} (\overline{\boldsymbol{\sigma}}) \;\; \text{with}\;\; \lambda >0  \quad \text{and} \quad \mathcal{D}^p(\overline{\boldsymbol{\sigma}}, \left\llbracket \boldsymbol{\epsilon}_p \right\rrbracket_{\bar{t}}) :=   \lambda \, |\overline{\boldsymbol{\sigma}}|. 
\end{equation}

The implementation of this example and the next two in \S \ref{int_hard_var} is achieved in \cite{Demoures2018c} through a discrete formulation of the theory developed in \S \ref{nonsmooth_elastoplasticity}. In Figures \ref{plasticity}, \ref{isot_plasticity}, \ref{kinemat_plasticity}, we plot the numerical approximations computed using a second-order multisymplectic variational integrator for elastoplastic problems. 

\begin{figure}[H] \centering 
 \includegraphics[width=1.51 in]{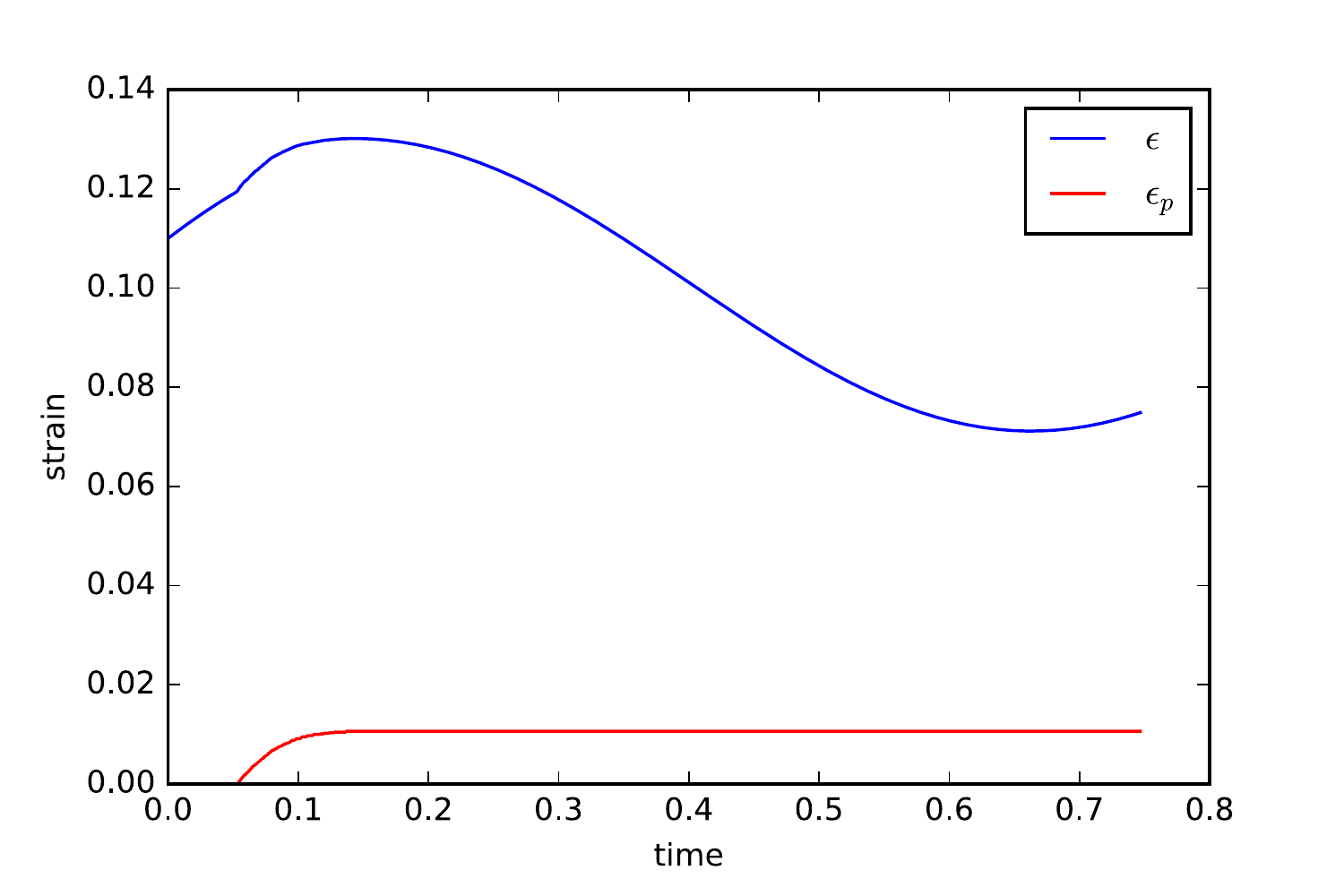} \vspace{-3pt}  \includegraphics[width=1.51 in]{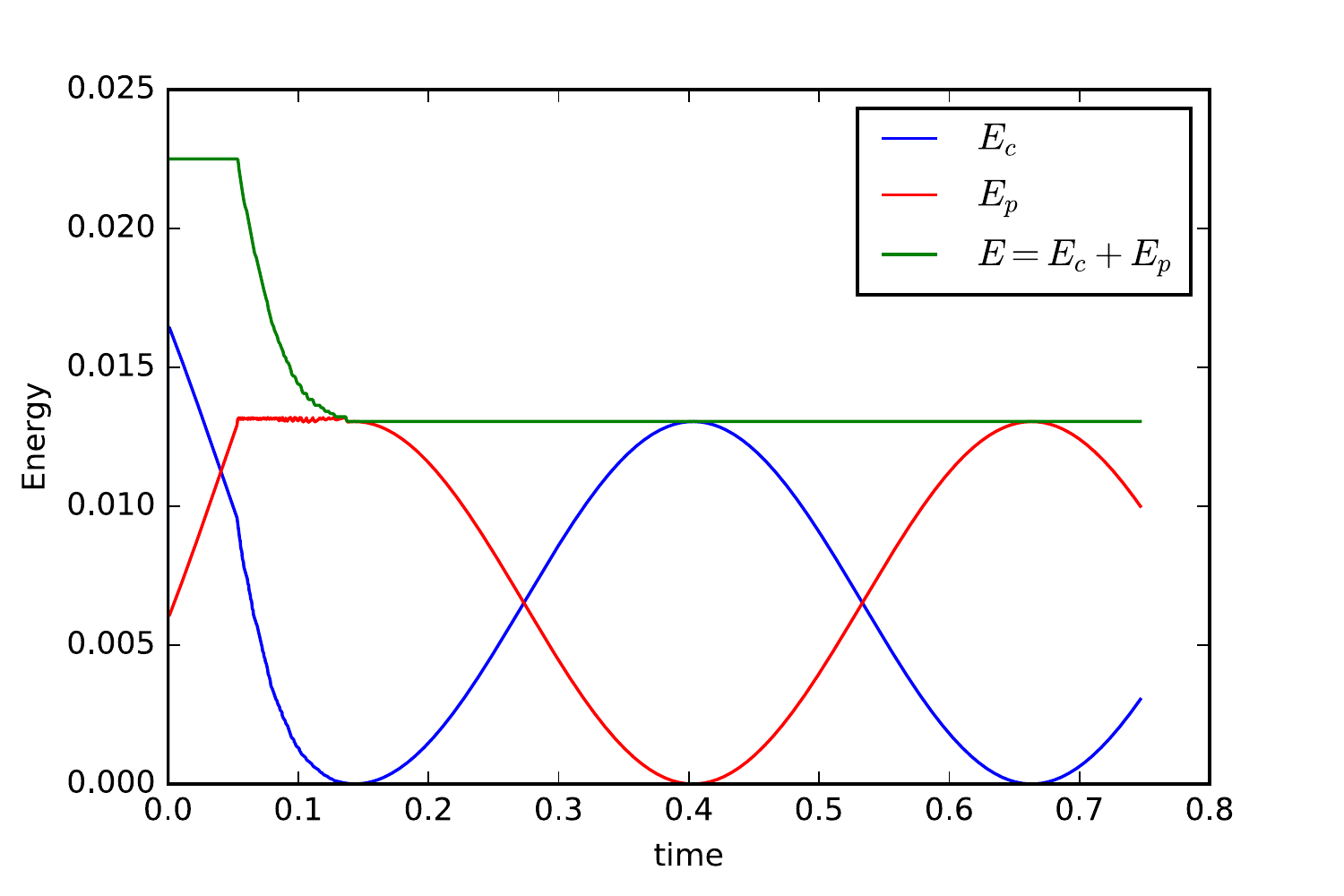} \vspace{-3pt} \includegraphics[width=1.51 in]{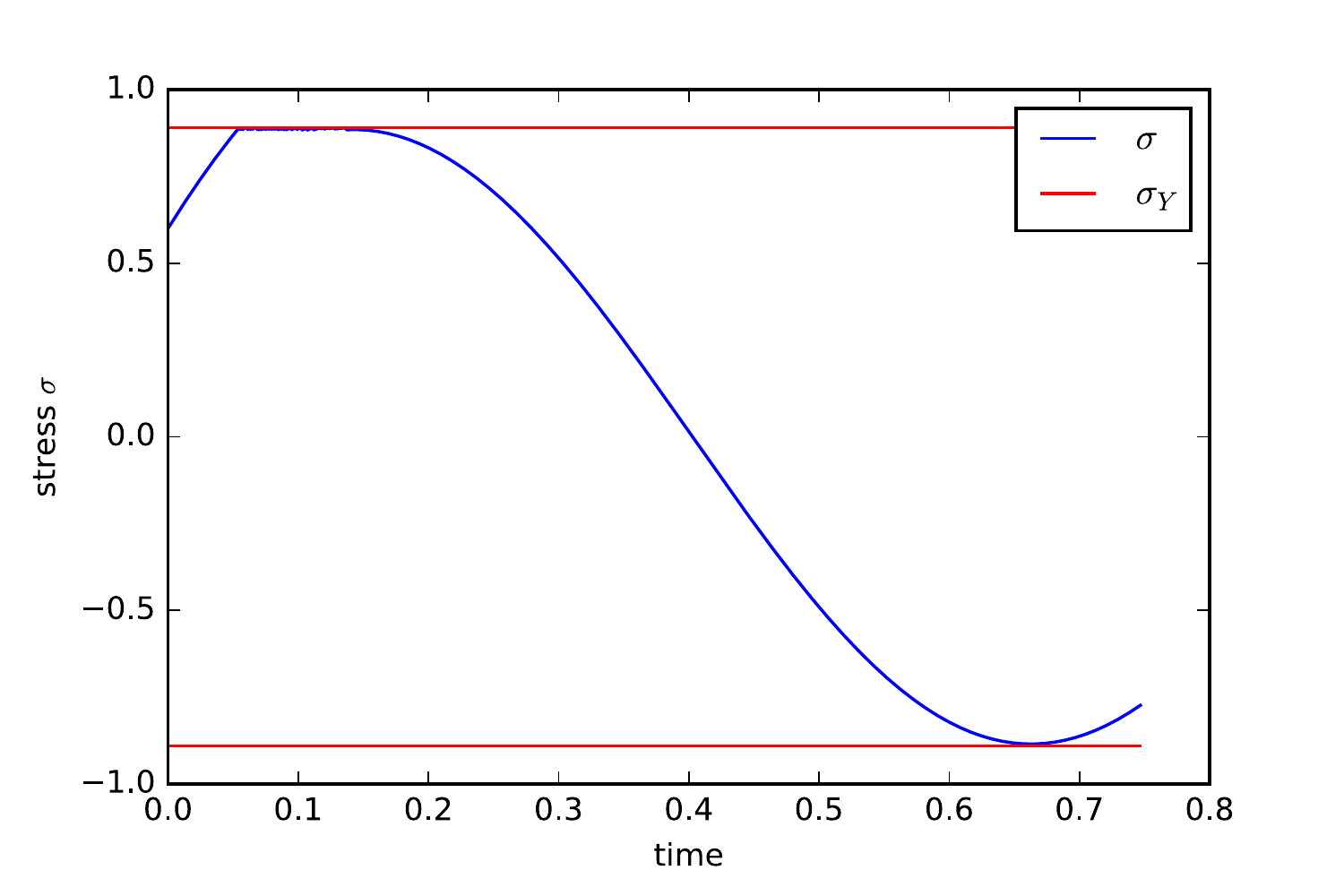} \vspace{-3pt}  \includegraphics[width=1.51 in]{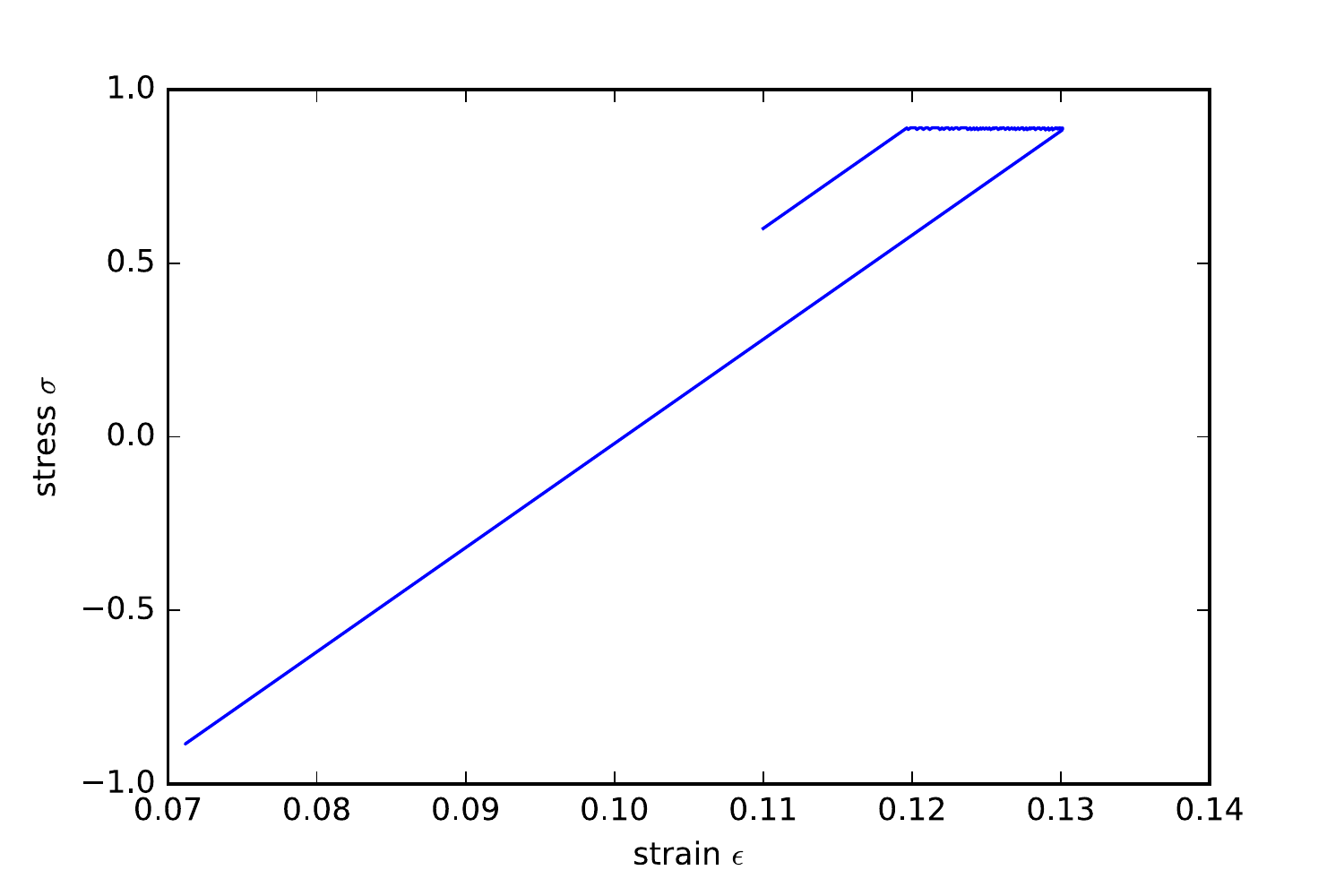} \vspace{-3pt}
 \caption{\footnotesize Tresca criterion. \textit{From left to right:} strain ($\boldsymbol{\epsilon}$, $\boldsymbol{\epsilon}_p$), total energy, stress, and stress/strain. $E=30$, $m=0.82$.} \label{plasticity} 
 \end{figure}
   
  In Figure \ref{plasticity} observe that the plastic strain $\boldsymbol{\epsilon}_p$ stops to evolve when $|\boldsymbol{\sigma}|  < \boldsymbol{\sigma}_Y$. Then after we go back to the elastic behavior and energy conservation. \quad $\lozenge$
 
\end{example}

\subsection{Internal hardening variables} \label{int_hard_var}

 Internal strain hardening variables $\boldsymbol{\xi}=(\xi_i,\xi_k)$ are often added to the plastic strain $\boldsymbol{\epsilon}_p$, where $\xi_i$, $\xi_k$ are respectively the isotropic and kinematic strain hardening variables. Then the potential energy is seen as the sum of the elastic store energy function \eqref{stress} plus the potential function $\mathcal{H}(\boldsymbol{\xi})$ for the hardening variables.

The yield criterion associated to pure plasticity \eqref{elas_plas_dom} can be modified in two ways. 
\noindent  \textbf{a)} \textit{Isotropic hardening:} The yield surface expands with increasing stress. Such that the yield criterion on $(\boldsymbol{\sigma}, \beta_i)$ is now defined as 
\begin{equation}\label{isot_yield_crit}
 f(\boldsymbol{\sigma}, \beta_i) := F(\boldsymbol{\sigma}) - k_{\boldsymbol{\sigma}_Y}(\beta_i) \leqslant 0 \quad \text{
 where} \quad \beta_i = - \partial \mathcal{H}(\boldsymbol{\xi})/\partial \xi_i.
 \end{equation}
 
\noindent \textbf{b)} \textit{Kinematic hardening:} The yield surface with the same shape is translated in stress space, with the following yield criterion 
\begin{equation}\label{kin_yield_crit}
f(\boldsymbol{\sigma},\beta_k) := F(\boldsymbol{\sigma}, \beta_k)  -  \boldsymbol{\sigma}_Y  \leqslant 0 \quad \text{
 where} \quad \beta_k = - \partial \mathcal{H}(\boldsymbol{\xi})/\partial \xi_k .
\end{equation} 

Note that the isotropic and kinematic hardening are often combined. But in the following we will consider the two cases separately.

\subsubsection{Isotropic hardening} \label{isotropic_hard}

We consider now the possibility of an expansion of the yield surface due to the increasing flow stress $k_{\boldsymbol{\sigma}_Y}(\beta_i)$. See, e.g., in \cite{Simo1998} the following elementary model $k_{\boldsymbol{\sigma}_Y}(\beta_i) = \boldsymbol{\sigma}_Y - \beta_i$ where $\boldsymbol{\sigma}_Y > 0$ is given constant.

\medskip

 The admissible set is now defined as
\begin{equation}\label{elas_plas_hard_dom}
\mathbb{E}_{(\sigma,q)}:=  \{ (\boldsymbol{\sigma},\beta_i) \in \mathbb{R}^2 \, | \, f(\boldsymbol{\sigma},\beta_i)\leq 0 \},
\end{equation}
where the yield criterion $f(\boldsymbol{\sigma},\beta_i)$ for isotropic hardening verifies \eqref{isot_yield_crit}. The Lagrangian is now defined by
\begin{equation}\label{Lagrangian_hard}
L(\boldsymbol{\epsilon}, \dot{\boldsymbol{\epsilon}}) = \frac{1}{2} m |\dot{\boldsymbol{\epsilon}}|^2 - \frac{1}{2} E |\boldsymbol{\epsilon}- \boldsymbol{\epsilon}_p|^2 - \mathcal{H}(\xi_i).
\end{equation}
 Through the derivative of the action $\mathfrak{S}^{ns}(\boldsymbol{\epsilon})$ outside of plastic behavior we get the Euler-Lagrange equations \eqref{EL}. 
 
 \medskip
 
 When the plastic phenomenon occurs at time $\overline{t}$ from the vertical jump condition we get the same as the one obtained with perfect plasticity, i.e., \eqref{vert_jump_pp}. While from the plasticity law \eqref{non_smooth_plast} we obtain the following rate of change of plastic strain and of isotropic hardening
 \begin{equation}\label{sigma_constraint_hard}
\left( \left\llbracket \boldsymbol{\epsilon}_p \right\rrbracket_{\bar{t}} ,  \left\llbracket \xi_i \right\rrbracket_{\bar{t}} \right) \in N_{\mathbb{E}_{(\sigma,q)}}(\overline{\boldsymbol{\sigma}} ,\overline{\beta}_i) = \lambda \nabla_{(\boldsymbol{\sigma}, \beta_i)}f( \overline{\boldsymbol{\sigma}} , \overline{\beta}_i), \;\; \text{for all} \;\; \lambda\geq 0.
 \end{equation}
 The rate of change of the total energy \eqref{rate_energ_0} becomes
\begin{equation}\label{rate_energ_isohard}
\begin{aligned}
 \left\llbracket - E_{\rm tot}  \right\rrbracket &  =  \left<(\overline{\boldsymbol{\sigma}},\overline{\beta}_i) , \lambda_i \nabla_{(\boldsymbol{\sigma}, \beta_i)}f( \overline{\boldsymbol{\sigma}} , \overline{\beta}_i) \right>  =:  \mathcal{D}_{ih}^p, \quad \text{with} \quad \lambda_i > 0,
\end{aligned}
\end{equation}
where $ \mathcal{D}_{ih}^p$ is denoted the \textit{isotropic hardening plastic dissipation function}. 

 \begin{example}\label{exemple_2}
 Let the 1D yield criterion corresponding to isotropic hardening
 \begin{equation}\label{non_lin_yield}
 f(\boldsymbol{\sigma}, \beta_i) := |\boldsymbol{\sigma}|  + \beta_i - \boldsymbol{\sigma}_Y  \leq 0.
 \end{equation}
 where $\boldsymbol{\sigma}_Y$ is constant. We specify the potential function for isotropic hardening variables
\[
\mathcal{H}(\xi_i) =  \frac{1}{2}K \xi_i^2 \quad \text{with} \;\;  K>0. 
\]
So we get
 \begin{equation*} %\label{relation_hard_bis}
\begin{bmatrix}  \boldsymbol{\sigma} \\ \beta_i \end{bmatrix}   = \begin{bmatrix} E(\boldsymbol{\epsilon}-\boldsymbol{\epsilon}_p) \\ -K \xi_i \end{bmatrix}  \;\; \text{and} \;\; \begin{bmatrix} \left\llbracket \boldsymbol{\epsilon}_p \right\rrbracket_{\bar{t}} \\  \left\llbracket \xi_i \right\rrbracket_{\bar{t}} \end{bmatrix} \in N_{\mathbb{E}_{(\sigma,q)}}(\overline{\boldsymbol{\sigma}}, \overline{\beta}_i)  =\lambda  \begin{bmatrix} \mathrm{sgn}(\overline{\boldsymbol{\sigma}}) \\ 1 \end{bmatrix}, \;\; \text{for all}\quad \lambda \geq 0. 
 \end{equation*}
 
 The numerical tests implemented through a discrete formulation of isotropic hardening give us the following results

\begin{figure}[H] \centering 
 \includegraphics[width=1.51 in]{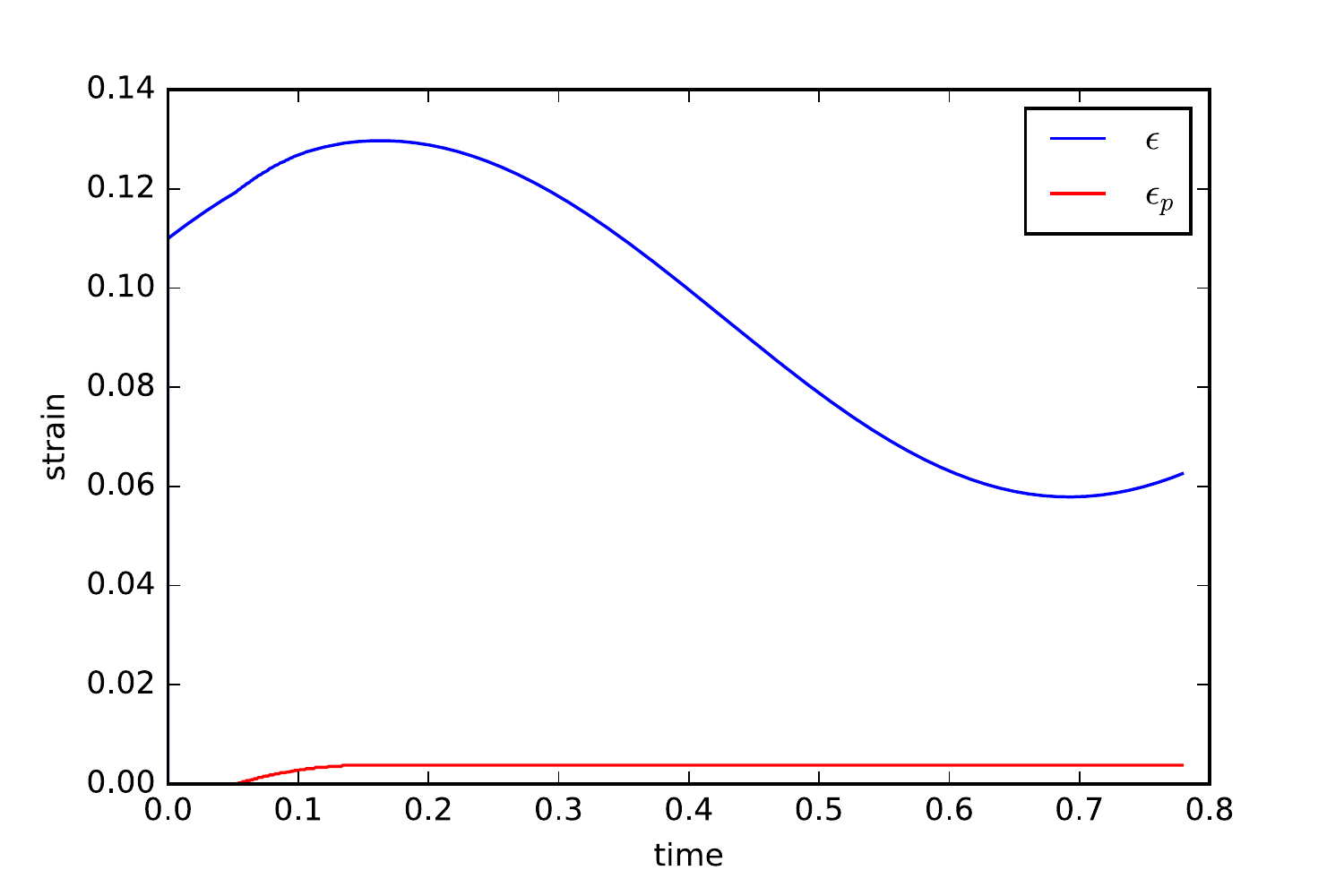} \vspace{-3pt}  \includegraphics[width=1.51 in]{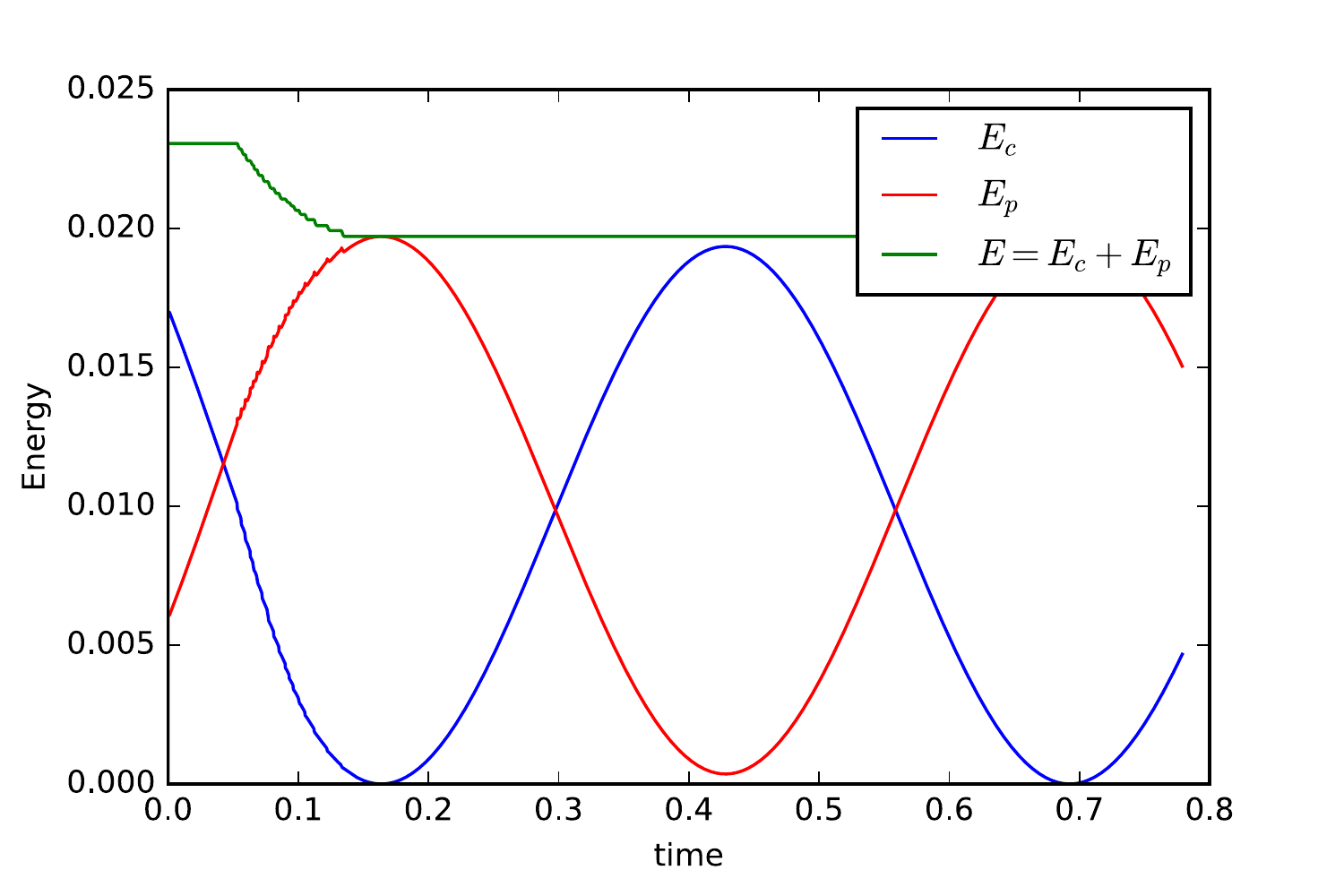} \vspace{-3pt} \includegraphics[width=1.51 in]{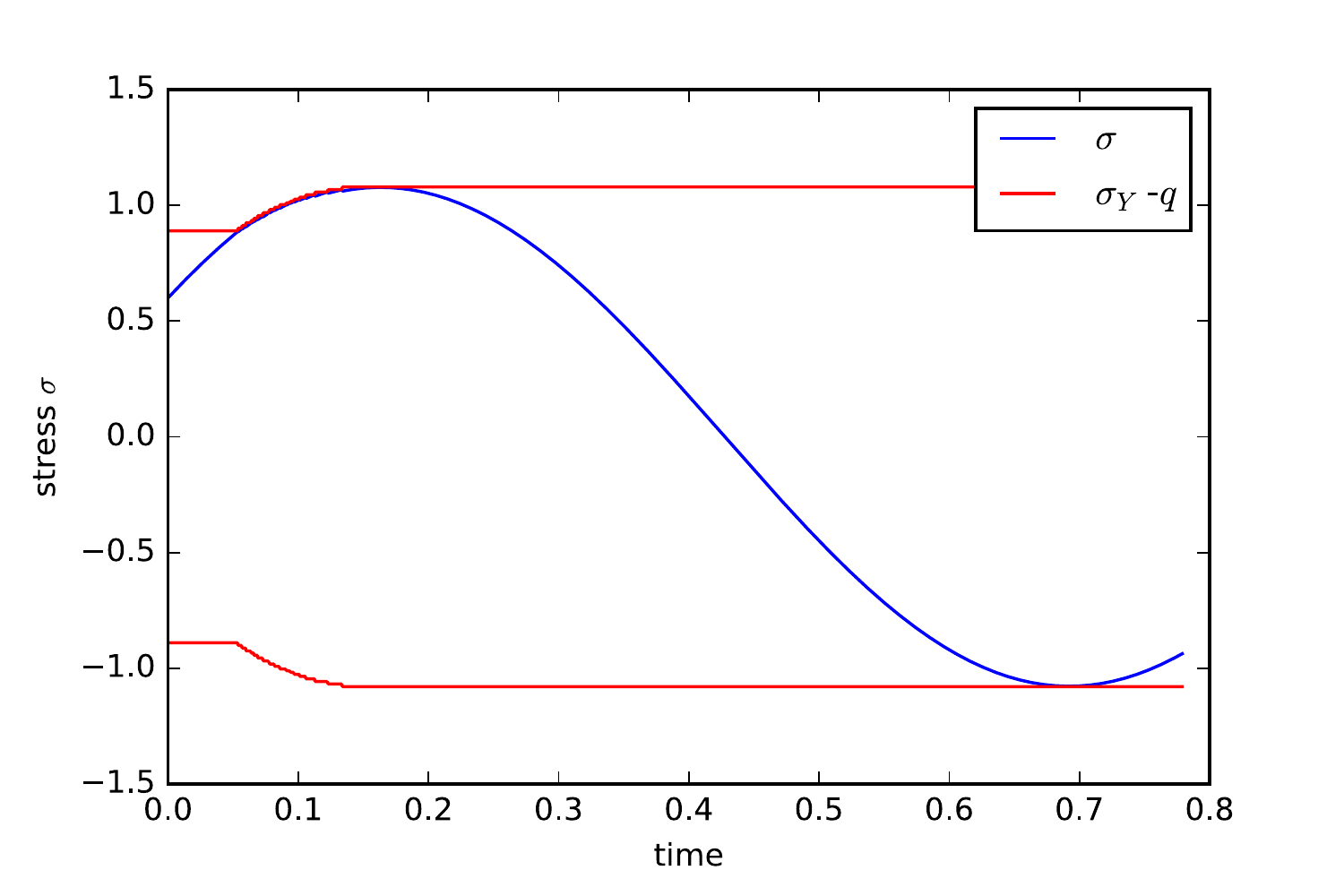} \vspace{-3pt}  \includegraphics[width=1.51 in]{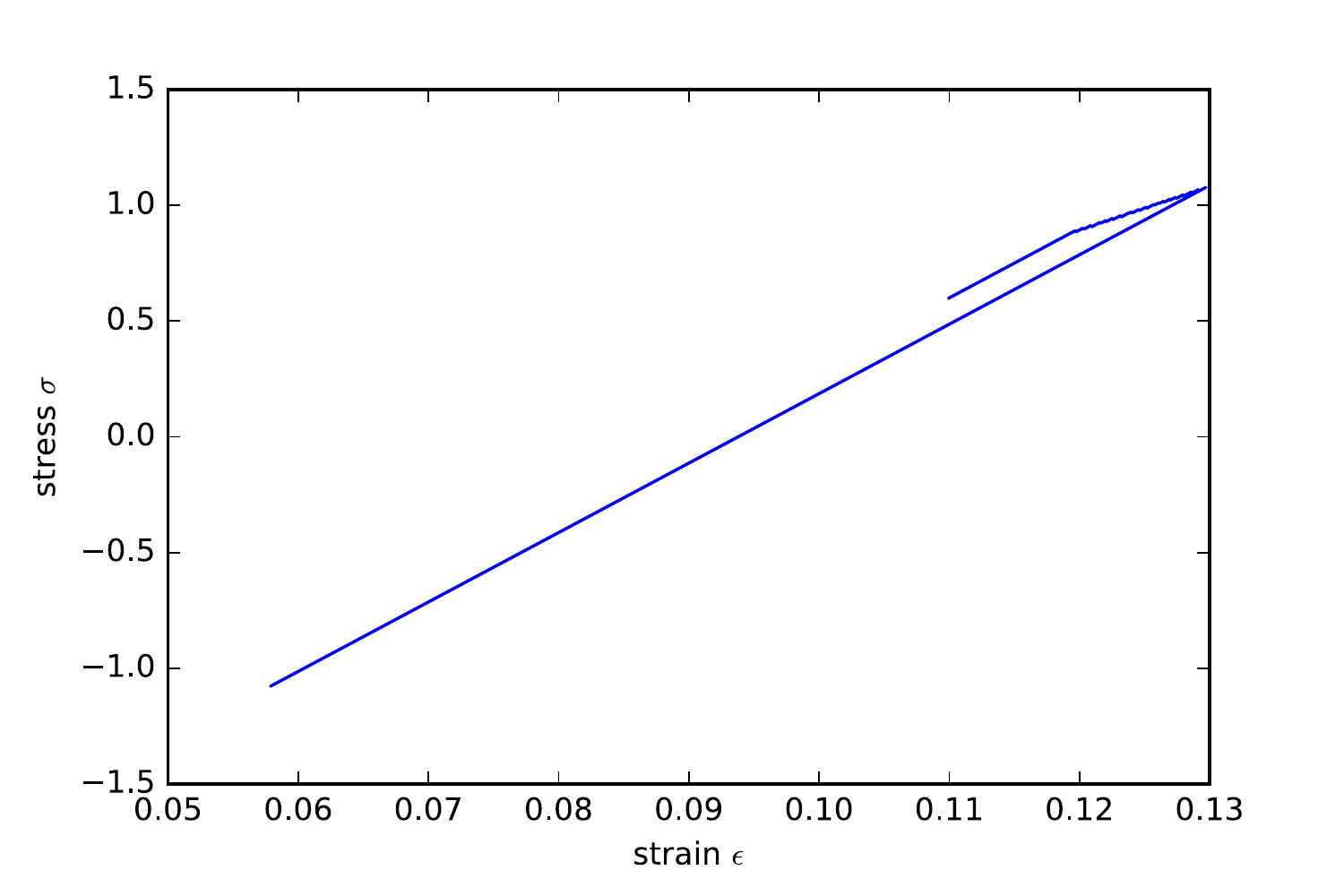} \vspace{-3pt}
 \caption{\footnotesize  \textit{From left to right:} strain ($\boldsymbol{\epsilon}$, $\boldsymbol{\epsilon}_p$), total energy, stress, and stress/strain. $E=30$, $K=50$, $m=0.85$. } \label{isot_plasticity} 
 \end{figure}
 When the stress satisfies $|\boldsymbol{\sigma}|=\boldsymbol{\sigma}_Y - \beta_i$, the plastic strain increases through small jumps, and stops as soon as $|\boldsymbol{\sigma}| < \boldsymbol{\sigma}_Y- \beta_i$. However note that the yield surface expands, due to $\boldsymbol{\sigma}_Y- \beta_i$ which increases. That is, loading after unloading will define a new instantaneous elastic limit and so forth. \quad $\lozenge$
 \end{example}
% We observe the Bauschinger effect which was discovered by Bauschinger [1886]

\subsubsection{kinematic hardening}\label{kinematic_hard}

The yield criterion $f(\boldsymbol{\sigma},\beta_k)$ described by \eqref{kin_yield_crit} exhibits \textit{kinematic hardening}. See, e.g., in \cite{Simo1998} the following elementary model $F(\boldsymbol{\sigma}, \beta_k) = \left|\boldsymbol{\sigma} - \beta_k \right|$. Where the yield surface translates in the direction of the plastic flow.

 When the plastic phenomenon occurs at time $\overline{t}$, the rate of change of plastic strain and kinematic hardening are given by 
 \begin{equation}\label{sigma_constraint_hard}
(\left\llbracket \boldsymbol{\epsilon}_p \right\rrbracket_{\bar{t}} , \left\llbracket \xi_k \right\rrbracket_{\bar{t}} ) \in N_{\mathbb{E}_{(\sigma,q)}}(\overline{\boldsymbol{\sigma}}, \overline{\beta}_k) = \lambda \nabla_{(\boldsymbol{\sigma}, \beta_k)}f( \overline{\boldsymbol{\sigma}}, \overline{\beta}_k), \;\; \text{for all} \;\; \lambda\geq 0.
 \end{equation}
 While the horizontal jump condition becomes
\begin{equation}\label{rate_energ_kinhard}
\begin{aligned}
 \left\llbracket - E_{\rm tot}  \right\rrbracket &  =  \left< (\overline{\boldsymbol{\sigma}}, \overline{\beta}_k) , \lambda_k \nabla_{(\boldsymbol{\sigma}, \beta_k)}f( \overline{\boldsymbol{\sigma}}, \overline{\beta}_k) \right> =:  \mathcal{D}_{kh}^p, \quad \text{with} \quad \lambda_k >0,
\end{aligned}
\end{equation}
where the \textit{kinematic hardening plastic dissipation function} is denoted by $ \mathcal{D}_{kh}^p$.

 \begin{example}\label{exemple_3}
  Let the 1D yield criterion corresponding to kinematic hardening
 \begin{equation}\label{kinem_hard}
 f(\boldsymbol{\sigma}, \beta_k) := |\boldsymbol{\sigma} - \beta_k|  - \boldsymbol{\sigma}_Y  \leq 0,
 \end{equation} 
 and the following potential function
\[
\mathcal{H}(\xi_k) =  \frac{1}{2}H(\xi_k)^2 \quad \text{with} \;\;  H>0. 
\]
We get
 \begin{equation*} %\label{relation_hard_bis}
\begin{bmatrix}  \boldsymbol{\sigma} \\ \beta_k \end{bmatrix}   = \begin{bmatrix} E(\boldsymbol{\epsilon}-\boldsymbol{\epsilon}_p) \\ -H \xi_k \end{bmatrix}   \;\; \text{and} \;\; \begin{bmatrix}\left\llbracket \boldsymbol{\epsilon}_p \right\rrbracket_{\bar{t}} \\ \left\llbracket \xi_k \right\rrbracket_{\bar{t}} \end{bmatrix} \in N_{\mathbb{E}_{(\sigma,\beta_k)}}(\overline{\boldsymbol{\sigma}}, \overline{\beta}_k)  =  \lambda   \begin{bmatrix} \mathrm{sgn}(\overline{\boldsymbol{\sigma}}- \overline{\beta}_k) \\ - \mathrm{sgn}(\overline{\boldsymbol{\sigma}}- \overline{\beta}_k) \end{bmatrix}, \;\; \forall \;\; \lambda \geq 0. 
 \end{equation*}
 
 The numerical tests implemented through a discrete formulation of kinematic hardening give us the following results
 \begin{figure}[H] \centering 
 \includegraphics[width=1.51 in]{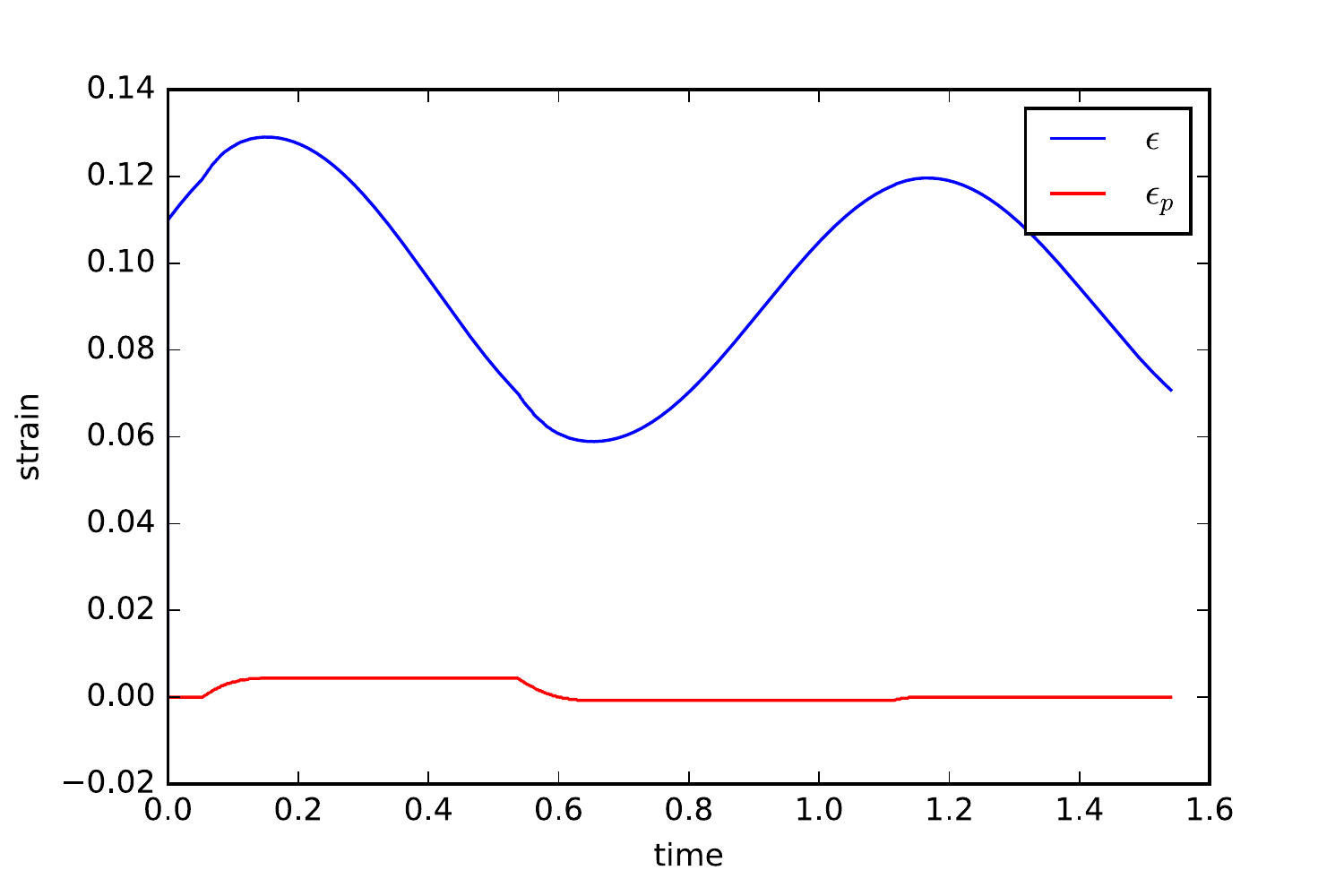} \vspace{-3pt}  \includegraphics[width=1.51 in]{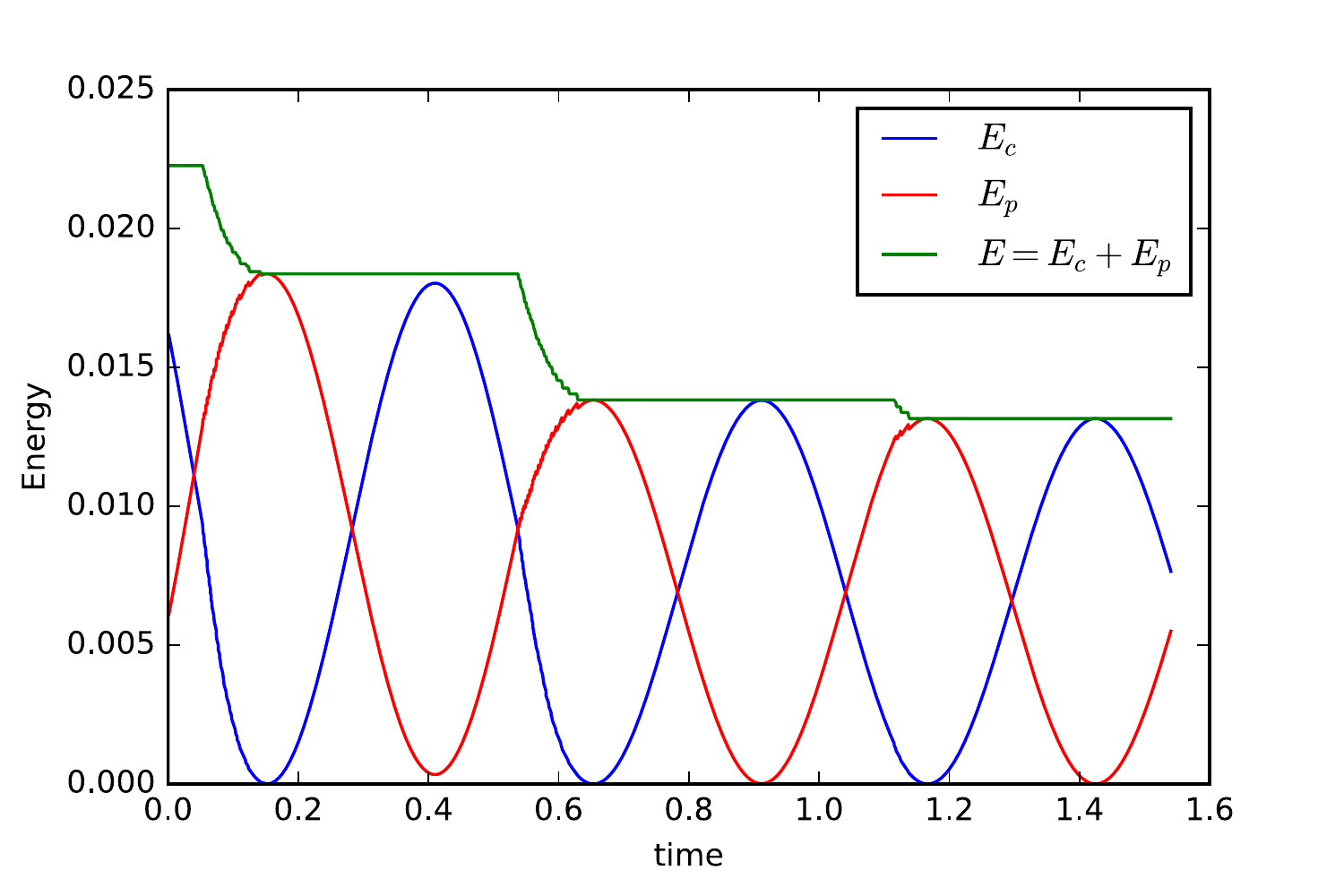} \vspace{-3pt} \includegraphics[width=1.51 in]{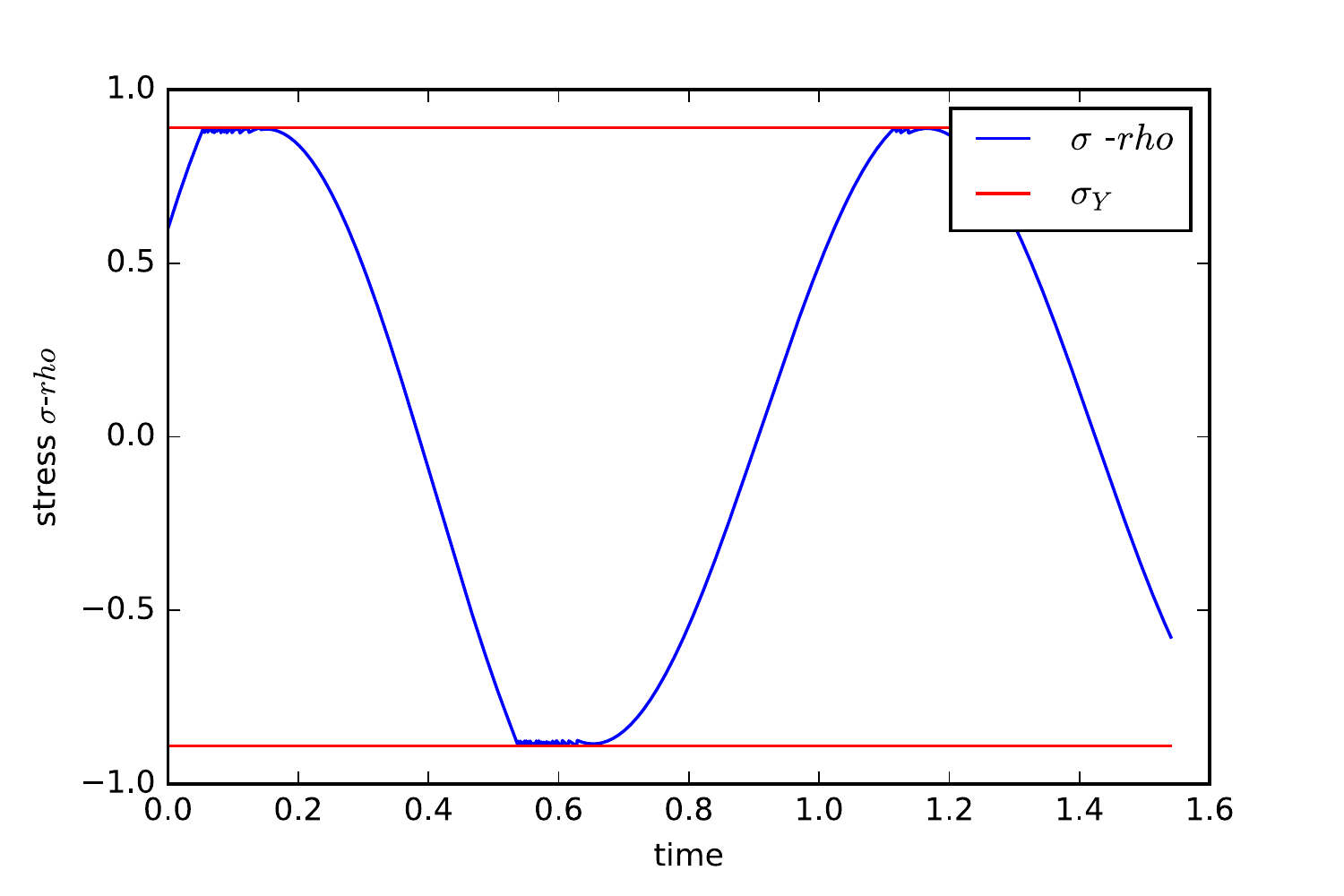} \vspace{-3pt}  \includegraphics[width=1.51 in]{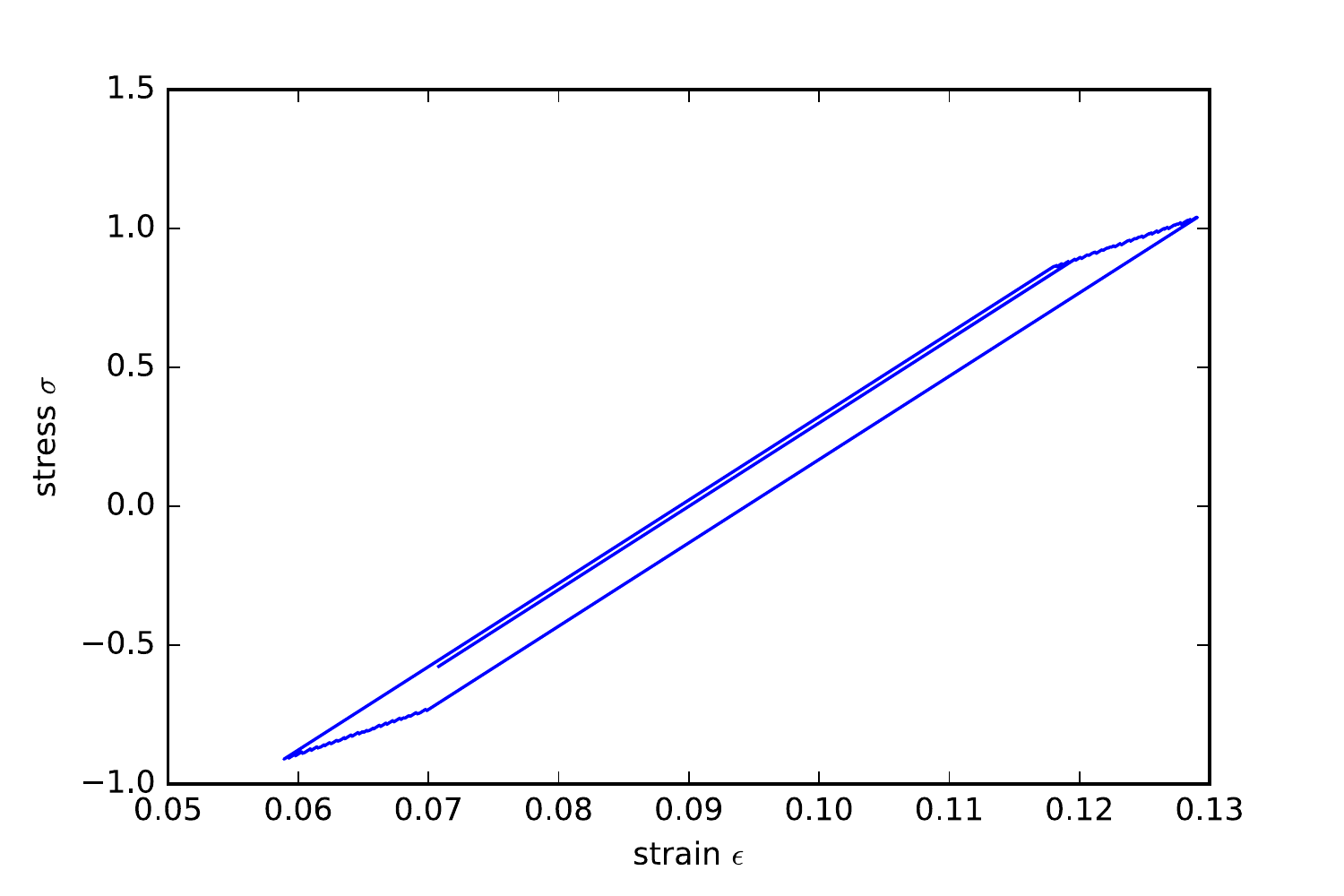} \vspace{-3pt}
 \caption{\footnotesize  \textit{From left to right:} strain ($\boldsymbol{\epsilon}$, $\boldsymbol{\epsilon}_p$), total energy, stress, and stress/strain. $E=30$, $H=35$, $m=0.81$. } \label{kinemat_plasticity} 
 \end{figure}
 Note that the yield surface retains the same shape but translates during the plastic strain, due to the internal variable $\beta_k$ in $|\boldsymbol{\sigma} - \beta_k|$. \quad $\lozenge$
 \end{example}

%%%%%%%%%%%%%%%%%%%%%%%%%%%%%%%%%%%%%%%%%%%%
%%%%%%%%%%%%%%%%%%%%%%%%%%%%%%%%%%%%%%%%%%%

\section{Rheological model for thermoplasticity}\label{rheo_nonsmooth_thermoelasto}

``First, suppose that we do irreversible work on an object by friction, generating a heat $Q$ on some object at temperature $T$. The entropy is increased by $Q/T$. The heat $Q$ is equal to the work $W$, and thus when we do a certain amount of work by friction against an object whose temperature is $T$, the entropy of the whole world increases by $W/T$.''(\cite{FeLeSa1963})

\medskip

The two main points which characterize the cristal plastic phenomenon are as follows: first, the heat is produced by the plastic strains; second, in the irreversible plastic change, the total entropy of the system always increases. 

\subsection{First and second laws of thermodynamics}\label{Cl_Du_dissip}

In the following we recall the second and first laws of thermodynamics. 

\paragraph{Second law.} The second law of thermodynamics was first put into words by \cite{Carnot1824}. It can be expressed as follows: 

In a isolated mechanical system which absorbes heat $Q_1$ at temperature $T_1$ and delivers heat $Q_2$ at temperature $T_2$, the relation between the two verifies
\begin{equation}\label{second_law}
\frac{Q_1}{T_1} = S = \frac{Q_2}{T_2},
\end{equation}
where $S$ denotes the entropy. 
However the second law of thermodynamics is expressed through different forms which depend from the perspective adopted. For example, the Clausius-Duhem local form of the ``second law of thermodynamics'' in a continuous body asserts that (see e.g., \cite{MaHu1994} \S 2.5)
\begin{equation} \label{Cl_Du}
 \gamma =   \dot{S} -  \frac{\rho r}{T} +  \frac{1}{T} \nabla \cdot \mathbf{q}  - \frac{1}{T^2} \nabla T \cdot \mathbf{q} \geq 0,
\end{equation}
where $\gamma $ is the rate of change of entropy production, $S$ is the entropy in the body, $r$ is the heat supply by unit of mass, $T$ is the temperature, $\mathbf{q}$ is the heat flux, and $\dot{S}$ is the rate of change of the total entropy. 

The dissipation $T \gamma$ is decomposed into the sum of the internal dissipation $\mathcal{D}_{int}$ under the Clausius-Plank form of the second law and the dissipation $\mathcal{D}_{cond}$ arising from heat conduction, see \cite{TrNo1965}(79.8, 79.9, 79.10), which are respectively 
\begin{equation}\label{D_int_D_cond}
\mathcal{D}_{int}:=  T \dot{S} - \rho r + \nabla \cdot \mathbf{q} \geq 0,   \quad \text{and} \quad \mathcal{D}_{cond}:= - \frac{1}{T} \nabla T \cdot \mathbf{q} \geq 0.
\end{equation} 
 In the expression \eqref{D_int_D_cond} note that we take into account the heat flux $\mathbf{q}$ which is function of the thermal conductivity of the material. 

\medskip

\paragraph{First law or balance of energy.} 
\cite{Clausius1850} and W. Rankine both stated the first law of thermodynamics which says that the rate of increase of the internal energy $\dot{E}_{int}$ of the body equals the rate of work done (the body forces and surface traction) plus the rate of increase of heat energy. 
\begin{equation}\label{first_law}
\dot{E}_{int} = \dot{W} + \dot{Q}.
\end{equation}

As a consequence, if we admit that the transformation is isothermal, then the work dissipated during plastic deformation is transformed into heat, i.e., the work lost by internal friction is equal to the heat produced.
\[
\dot{Q}_{plastic} = - \dot{W}_{friction} \quad \text{at temperature}\;\; T.
\] 

\medskip

We recall from \cite{MaHu1994} \S 2.3 the form it takes in a continuous body 
\begin{equation}\label{en_bal_eq}
\dot{E}_{int} = \boldsymbol{\sigma}: \mathbf{d} + \rho r - \nabla \cdot \mathbf{q},
\end{equation}
where $\boldsymbol{\sigma}$ is the Cauchy stress tensor and $ \mathbf{d}$ is the rate of change of strain tensor.

\medskip

By combination of the Clausius-Duhem inequality $\eqref{D_int_D_cond}_1$ and the energy balance \eqref{en_bal_eq} we can express the \textit{internal dissipation} $\mathcal{D}_{int}$ in the solids as the thermal power plus the mechanical power minus the time rate of change of the internal energy, i.e.,
\begin{equation} \label{int_dissip}
\mathcal{D}_{int} = T \dot{S} + \boldsymbol{\sigma}: \mathbf{d} - \dot{E}_{int} \geq 0.
\end{equation}

\subsection{Thermoplasticity: Perfectly plastic case} \label{thermo_perf_plast}

\paragraph{Thermomechanics of plasticity.} We introduce temperature $T$ as a new variable in the 1D elastoplastic model described in section \S\ref{elas_plast_mod}. From thermodynamics principle the stored energy \eqref{stress} becomes a function of the elastic strain and of the elastic entropy. Which is now denoted $\mathscr{W}(\boldsymbol{\epsilon}_e,S_e)$. Then, we recall the following \textit{local state axioms}; see e.g., \cite{Maugin1992}
\begin{equation} \label{sigma_theta}
\boldsymbol{\sigma}= \partial_{\boldsymbol{\epsilon}_e}\mathscr{W}(\boldsymbol{\epsilon}_e,S_e), \qquad T= \partial_{S_e}\mathscr{W}(\boldsymbol{\epsilon}_e,S_e).
\end{equation}

We admit that the total entropy $S$ is the sum
 of the entropy $S_e$ due to the elasticity and the entropy $S_p$ due to the plasticity, i.e.,
 \begin{equation}\label{entropy}
S=S_e + S_p.
\end{equation}

The \textit{Helmholtz free energy} $\Psi(\boldsymbol{\epsilon}_e, T)$ is defined from $\mathscr{W}$ by performing the change of variable $S_e \rightarrow T$ through the Legendre-Fenchel transform, see e.g. \cite{Ottinger2005}
\begin{equation}\label{free_energy_thermo}
\Psi(\boldsymbol{\epsilon}_e, T) = \mathscr{W}(\boldsymbol{\epsilon}_e,S_e) - T S_e .
\end{equation}
Hence, the local state axioms \eqref{sigma_theta} are now expressed as follows
\begin{equation} \label{free_entropy}
\boldsymbol{\sigma}= \partial_{\boldsymbol{\epsilon}_e}\Psi(\boldsymbol{\epsilon}_e,T), \qquad  S_e = - \partial_T \Psi(\boldsymbol{\epsilon}_e,T).
\end{equation}

%\begin{example}: As described in \cite{Simo1998}, let the following thermomechanical extension of the purely mechanical model characterized by the stored energy function $W(\boldsymbol{\epsilon}_e)$, such that the free energy $\Psi$ takes the following form
%\begin{equation}\label{free_energy_ex}
%\begin{aligned} \Psi(\boldsymbol{\epsilon}_e,T) &= g(T) - (T-T_0)M(\boldsymbol{\epsilon}_e) +W(\boldsymbol{\epsilon}_e), \\ \text{with} \;\; g(T) & = c_0\left[(T -T_0) - T \, \mathrm{log}(T/ T_0) \right],\end{aligned} \end{equation} where $T_0$ is a reference temperature, $c_0$ is the heat capacity, and $(T-T_0)M(\boldsymbol{\epsilon}_e)$ is a thermo-mechanical coupling. From \eqref{free_entropy} we get the elastic entropy \[ S_e = c_0 \, \mathrm{log}(T/T_0) + M(\boldsymbol{\epsilon}_e). \qquad \lozenge \]
%\end{example}

\paragraph{Thermoelastic regime.} The mechanical system is described by its Lagrangian composed of the kinetic energy minus the Helmholtz free energy 
\begin{equation}\label{thermo_Lagrangian}
\mathcal{L}_T(\boldsymbol{\epsilon}, \dot{\boldsymbol{\epsilon}}) = \frac{1}{2} m |\dot{\boldsymbol{\epsilon}}|^2 - \Psi(\boldsymbol{\epsilon} - \boldsymbol{\epsilon}_p,T).
\end{equation}
Through the derivative of the action $\mathfrak{S}^{ns}(\boldsymbol{\epsilon})= \int \mathcal{L}_T(\boldsymbol{\epsilon}, \dot{\boldsymbol{\epsilon}})$ with respect to $\boldsymbol{\epsilon}$
we get the Euler-Lagrange equation 
\begin{equation}\label{thermo_EL}
m\ddot{\boldsymbol{\epsilon}} + \partial_{\boldsymbol{\epsilon}}\Psi(\boldsymbol{\epsilon} - \boldsymbol{\epsilon}_p,T)=0. 
\end{equation}

\begin{remark}{\rm
Note that we could define the Lagrangian \eqref{thermo_Lagrangian} by taking into account the internal energy $\mathscr{W}(\boldsymbol{\epsilon}_e,S_e)$ instead of the free energy \eqref{free_energy_thermo}, with a similar result. \quad $\square$ }
\end{remark}

\paragraph{Internal dissipation.} The general constitutive equations \eqref{int_dissip} provides the internal plastic dissipation $\mathcal{D}_{int}$, i.e.,  
\begin{equation} \label{int_dissip_0}
\begin{aligned}
\mathcal{D}_{int} & = T( \dot{S}_e + \left\llbracket S_p \right\rrbracket_t )  +   \left< \boldsymbol{\sigma}, \dot{\boldsymbol{\epsilon}}_e +\left\llbracket \boldsymbol{\epsilon}_p \right\rrbracket_t \right> - \dot{\mathscr{W}}(\boldsymbol{\epsilon}_e,S_e ) %\nonumber
\\
&  \stackrel{\eqref{sigma_theta}}=  T \left\llbracket S_p \right\rrbracket_t + \left<  \boldsymbol{\sigma}, \left\llbracket \boldsymbol{\epsilon}_p \right\rrbracket_t \right> \geq 0,   
\end{aligned}
\end{equation}
where the plastic strain $\left\llbracket \boldsymbol{\epsilon}_p \right\rrbracket_t$ evolves by jumps, and therefore the rate of change of plastic entropy $\left\llbracket S_p \right\rrbracket_t$ evolves also by jumps. 

We deduce the decomposition of the internal dissipation \eqref{int_dissip_0} into \textit{mechanical dissipation} $\mathcal{D}_{mech}$ and \textit{thermic dissipation} $\mathcal{D}_{ther}$, which are respectively
\begin{equation}\label{int_dissip_decomp}
\mathcal{D}_{mech}(\boldsymbol{\sigma}, \left\llbracket \boldsymbol{\epsilon}_p \right\rrbracket_t) := \left<  \boldsymbol{\sigma}, \left\llbracket \boldsymbol{\epsilon}_p \right\rrbracket_t \right> \quad \text{and} \quad \mathcal{D}_{ther}(T ,\left\llbracket S_p \right\rrbracket_t) := T \left\llbracket S_p \right\rrbracket_t.
\end{equation}

\paragraph{Thermoelastoplastic domain.} In the context of thermoplasticity the internal plastic dissipation \eqref{int_dissip_0} verifies the maximum dissipation principle; see \cite{Lubliner1984}, \cite{Simo1998} \S57. 

The elastoplastic domain \eqref{elas_plas_dom} is modified through the introduction of the temperature $T$. Therefore, the thermoelastoplastic domain is defined as follows
\begin{equation}\label{elas_plas_dom_theta}
\mathbb{E}_{\sigma, T}:= \left\{(\boldsymbol{\sigma}, T)\, | \, f(\boldsymbol{\sigma}, T) \leq 0 \right\},
\end{equation}
where $\mathrm{Int}(\mathbb{E}_{\sigma, T})$ and $\partial \mathbb{E}_{\sigma, T}$ define respectively the thermoelastic and the thermoplastic domains.

Given the maximum dissipation principle the problem we have to solve is to minimize $-\mathcal{D}_{mech}$ and $- \mathcal{D}_{ther}$ under the constraint $f(\boldsymbol{\sigma}, T) \leq 0$. From the results recalled in \S \ref{Moreau}, by extension, we deduce the following expressions
\begin{equation} \label{Normal_cone_theta}
\begin{aligned}
& \begin{bmatrix} \left\llbracket \boldsymbol{\epsilon}_p \right\rrbracket_t \\ \left\llbracket S_p \right\rrbracket_t \end{bmatrix}  \in N_{\mathbb{E}_{\sigma, T}}( \boldsymbol{\sigma}, T)= \lambda \begin{bmatrix} \partial_{\boldsymbol{\sigma}}  f(\boldsymbol{\sigma}, T) \\ \partial_T  f(\boldsymbol{\sigma}, T) \end{bmatrix} , \;\; \text{for all} \;\;  \lambda \geq 0, 
\\
&\text{or equivalently} \quad \begin{bmatrix} \boldsymbol{\sigma} \\ T \end{bmatrix}  \in N_{\mathbb{E}_{\sigma, T}}^*( \left\llbracket \boldsymbol{\epsilon}_p \right\rrbracket_t , \left\llbracket S_p \right\rrbracket_t ).
\end{aligned}
\end{equation}
From \eqref{Normal_cone_theta} we get $(\left\llbracket \boldsymbol{\epsilon}_p \right\rrbracket_t , \left\llbracket S_p \right\rrbracket_t )$.

\paragraph{Elastic entropy.} The expression \eqref{int_dissip_0} issued from \eqref{int_dissip} together with the constitutive equations $\eqref{D_int_D_cond}_1$ yields the following relation
\begin{equation} \label{Import_relat}
T \dot{S}_e = \left<\boldsymbol{\sigma}, \left\llbracket \boldsymbol{\epsilon}_p \right\rrbracket_t \right> - \nabla \cdot \mathbf{q} .
\end{equation}

By assumption we admit that there is no heat flux nor heat supply in the 1D elastoplastic model described in section \S\ref{elas_plast_mod}, i.e., $\nabla \cdot  \mathbf{q}=0$. So we get the following rate of change of elastic entropy
\begin{equation}\label{T_dotSe}
T \dot{S}_e=  \left< \boldsymbol{\sigma}, \left\llbracket \boldsymbol{\epsilon}_p \right\rrbracket_t \right> \quad \Leftrightarrow \quad  \dot{S}_e=  \frac{1}{T} \left< \boldsymbol{\sigma}, \left\llbracket \boldsymbol{\epsilon}_p \right\rrbracket_t \right>.
\end{equation}

The mechanical dissipation $\eqref{int_dissip_decomp}_1$ is described by jumps of energy, see \S \ref{Moreau}. As a consequence, from \eqref{T_dotSe}, we can conclude that the rate of change of the elastic entropy $\dot{S}_e$ evolve by jumps.
From now on the rate of change of elastic entropy is denoted 
\begin{equation} \label{rate_el_entropy}
\left\llbracket S_e \right\rrbracket_t: = \frac{1}{T} \left< \boldsymbol{\sigma},\left\llbracket \boldsymbol{\epsilon}_p \right\rrbracket_t \right> .
\end{equation} 

In a logical way, outside of plastic phenomenon, the rate of change of elastic entropy \eqref{rate_el_entropy} becomes $\left\llbracket S_e \right\rrbracket_t =0$. 

\paragraph{Total energy and entropy production.} The total energy is composed of the kinetic energy $\frac{1}{2} m |\dot{\boldsymbol{\epsilon}}|^2$, associated to the rheological model defined in \S \ref{elas_plast_mod}, plus the internal energy $\mathscr{W}(\boldsymbol{\epsilon}_e,S_e)$. The rate of change of the total energy $\left\llbracket E_{\rm tot} \right\rrbracket_t$, can be derive from \eqref{int_dissip_0}, i.e.,
\begin{equation}\label{rate_int_energ_ther}
\begin{aligned}
\left\llbracket E_{\rm tot} \right\rrbracket_t & \stackrel{\eqref{sigma_theta}}  = \dot{\boldsymbol{\epsilon}} m\ddot{\boldsymbol{\epsilon}} +\left< \overline{\boldsymbol{\sigma}} , \dot{\boldsymbol{\epsilon}}_e\right>+ T\, \left\llbracket S_e \right\rrbracket_t \\
&   \stackrel{ \eqref{total_displ} \eqref{thermo_EL}} = -  \left< \overline{\boldsymbol{\sigma}}, \left\llbracket \boldsymbol{\epsilon}_p \right\rrbracket_t \right>+  T \, \left\llbracket S_e \right\rrbracket_t    \stackrel{\eqref{rate_el_entropy}}= 0,
\end{aligned}
\end{equation}
where we used the relation $\dot{\boldsymbol{\epsilon}}=\dot{\boldsymbol{\epsilon}}_e + \left\llbracket \boldsymbol{\epsilon}_p \right\rrbracket_t$.
Hence the total energy is conserved. Finally the rate of change of the entropy production \eqref{Cl_Du} has the following expression
\begin{equation}\label{entrop_prod}
\begin{aligned}
\gamma & = \left\llbracket S_e \right\rrbracket_t  + \left\llbracket S_p \right\rrbracket_t  = \frac{1}{T} \left< \overline{\boldsymbol{\sigma}}, \left\llbracket \boldsymbol{\epsilon}_p \right\rrbracket_t \right> +  \left\llbracket S_p \right\rrbracket_t 
\\
& \stackrel{\eqref{Normal_cone_theta}} = \frac{\lambda}{T} \left< \overline{\boldsymbol{\sigma}} , \partial_{\boldsymbol{\sigma}}  f(\overline{\boldsymbol{\sigma}}, T) \right> +  \lambda \, \partial_T  f(\overline{\boldsymbol{\sigma}}, T),
\end{aligned}
\end{equation}
where $\gamma$ evolves by jumps.

\begin{example} : 
Let the following yield criterion issued from Tresca criterion \eqref{Tresca_law}
\[
f(\boldsymbol{\sigma}, T) = |\boldsymbol{\sigma}| - \boldsymbol{\sigma}_Y(T) \leqslant 0,
\]

During the plastic phenomena, when $f(\overline{\boldsymbol{\sigma}}, T)=0$, we get the mechanical dissipation and the thermic dissipation from \eqref{Normal_cone_theta} :
\[
\mathcal{D}_{mech}(\overline{\boldsymbol{\sigma}}, \left\llbracket \boldsymbol{\epsilon}_p \right\rrbracket_t)= \lambda \, | \overline{\boldsymbol{\sigma}} |,  \quad \text{and} \quad \mathcal{D}_{ther} = - \lambda T \frac{d}{d T} \boldsymbol{\sigma}_Y(T).  
\]
In addition we obtain 
\[
\gamma = \frac{\lambda}{T} \, | \overline{\boldsymbol{\sigma}} | - \lambda  \frac{d}{d T} \boldsymbol{\sigma}_Y(T) \geq 0. \qquad \lozenge
\]
\end{example}

\subsection{Hardening laws in thermoplasticity}

\subsubsection{Isotropic hardening } \label{thermo_isot_hard}

\paragraph{Thermomechanics of plasticity.} Let the rheological model with isotropic hardening law, as described in \S \ref{isotropic_hard}, where the internal energy is given by \eqref{Lagrangian_hard}. If temperature and entropy are introduced in it, the internal energy becomes 
\begin{equation}\label{int_energy_hard_temp}
 \mathscr{W}_{ih}(\boldsymbol{\epsilon}_e, \xi_i, S_e) = \mathscr{W}(\boldsymbol{\epsilon}_e,S_e) + \mathcal{H}(\xi_i),
\end{equation}
whereas the relationships \eqref{sigma_theta} are transformed into
\begin{equation} \label{sigma_theta_q}
\boldsymbol{\sigma}= \partial_{\boldsymbol{\epsilon}_e} \mathscr{W}_{ih}(\boldsymbol{\epsilon}_e, \xi_i, S_e), \;\;  \beta_i =- \partial_{\xi_i}  \mathscr{W}_{ih}(\boldsymbol{\epsilon}_e, \xi_i, S_e), \;\; T= \partial_{S_e} \mathscr{W}_{ih}(\boldsymbol{\epsilon}_e, \xi_i, S_e).
\end{equation}
Then, the Helmotz free energy is now defined as follows
\begin{equation} \label{free_energy_thermo_hard}
\Psi_{ih} (\boldsymbol{\epsilon}_e, \xi_i, T) =  \mathscr{W}_{ih}(\boldsymbol{\epsilon}_e, \xi_i, S_e) -TS_e, \quad \text{with}
\end{equation} 
\begin{equation} \label{free_entropy_q}
\boldsymbol{\sigma}= \partial_{\boldsymbol{\epsilon}_e}\Psi_{ih}(\boldsymbol{\epsilon}_e, \xi_i, T), \quad  \beta_i= - \partial_{\xi_i}\Psi_{ih}(\boldsymbol{\epsilon}_e, \xi_i, T), \quad S_e = - \partial_T \Psi_{ih}(\boldsymbol{\epsilon}_e, \xi_i, T).
\end{equation}

\paragraph{Internal dissipation.} For the isotropic hardening law we get the internal dissipation $\mathcal{D}_{int}$ from the general constitutive equations \eqref{int_dissip} and \eqref{sigma_theta_q}, i.e.,
\begin{equation} \label{int_dissip_hard}
\mathcal{D}_{int} = \left<  \boldsymbol{\sigma}, \left\llbracket \boldsymbol{\epsilon}_p \right\rrbracket_t \right> + \left< \beta_i, \left\llbracket \xi_i \right\rrbracket \right>  + T \left\llbracket S_p \right\rrbracket_t  = \mathcal{D}_{mech} + \mathcal{D}_{ther} \geq 0,  
\end{equation}
which verify the maximum dissipation principle. The expression of the rate of change of the elastic entropy is obtained by \eqref{int_dissip_hard} together with the constitutive equations $\eqref{D_int_D_cond}_1$. We get 
\begin{equation}\label{dotSe_isohar}
\left\llbracket S_e \right\rrbracket_t = \frac{1}{T} \big( \left<  \boldsymbol{\sigma}, \left\llbracket \boldsymbol{\epsilon}_p \right\rrbracket_t  \right> + \left< \beta_i, \left\llbracket \xi_i \right\rrbracket_t \right>\big).
\end{equation}

\paragraph{Thermoelastoplastic domain.} It is now defined as follows
\begin{equation}\label{elas_plas_dom_theta_hard}
\mathbb{E}_{ih}:= \left\{(\boldsymbol{\sigma}, \beta_i,T)\, | \, f_{ih}(\boldsymbol{\sigma}, \beta_i,T) \leq 0 \right\}.
\end{equation}
The problem to solve becomes: to minimize $- \mathcal{D}_{mech}$ and $- \mathcal{D}_{ther}$ under the constraint $f_{ih}(\boldsymbol{\sigma}, \beta_i,T)\leq 0$. We obtain
\begin{equation} \label{Normal_cone_theta_hard}
\begin{aligned}
& \begin{bmatrix} \left\llbracket \boldsymbol{\epsilon}_p \right\rrbracket_t \\ \left\llbracket \xi_i \right\rrbracket_t  \\ \left\llbracket S_p \right\rrbracket_t  \end{bmatrix}  \in N_{\mathbb{E}_h}( \boldsymbol{\sigma}, \beta_i, T)= \lambda \begin{bmatrix} \partial_{\boldsymbol{\sigma}}  f_{ih}(\boldsymbol{\sigma}, \beta_i , T) \\ \partial_{\beta_i}  f_{ih}(\boldsymbol{\sigma}, \beta_i , T) \\ \partial_T  f_{ih}(\boldsymbol{\sigma}, \beta_i , T) \end{bmatrix} , \;\; \text{for all} \;\;  \lambda \geq 0, 
\\
& \text{or equivalently} \quad  \begin{bmatrix} \boldsymbol{\sigma} & \beta_i & T \end{bmatrix}^\mathrm{T}  \in N_{\mathbb{E}_h}^*(\left\llbracket \boldsymbol{\epsilon}_p \right\rrbracket_t , \left\llbracket \xi_i \right\rrbracket_t , \left\llbracket S_p \right\rrbracket_t ).
\end{aligned}
\end{equation}

\paragraph{Variational formulation.} The Lagrangian composed of the kinetic energy minus the free energy
\begin{equation}\label{thermo_Lagrangian_2}
\mathcal{L}_T(\boldsymbol{\epsilon}, \dot{\boldsymbol{\epsilon}}) = \frac{1}{2} m |\dot{\boldsymbol{\epsilon}}|^2 - \Psi_{ih}(\boldsymbol{\epsilon} - \boldsymbol{\epsilon}_p, \xi_i, T).
\end{equation}
The CEL equation is obtained through the derivative of the action with respect of $\boldsymbol{\epsilon}$. The resulting expression is 
\begin{equation}\label{thermo_ELisot}
m\ddot{\boldsymbol{\epsilon}} + \partial_{\boldsymbol{\epsilon}}\Psi_{ih}(\boldsymbol{\epsilon} - \boldsymbol{\epsilon}_p, \xi_i, T)=0. 
\end{equation}

\paragraph{Entropy production.} At temperature $T$ the total energy is conserved and the rate of change of entropy production $ \gamma$ becomes
 \begin{equation}\label{entrop_prod}
\begin{aligned}
\gamma & = \left\llbracket S_e \right\rrbracket_t  + \left\llbracket S_p \right\rrbracket_t = \frac{1}{T} \Big( \left<  \overline{\boldsymbol{\sigma}}, \left\llbracket \boldsymbol{\epsilon}_p \right\rrbracket_t \right> + \left< \overline{\beta}_i, \left\llbracket \xi_i \right\rrbracket_t \right>\Big) + \left\llbracket S_p \right\rrbracket_t
\\
& \stackrel{\eqref{Normal_cone_theta_hard}} = \frac{\lambda}{T} \Big( \left< \overline{\boldsymbol{\sigma}}, \partial_{\boldsymbol{\sigma}}  f_{ih}(\overline{\boldsymbol{\sigma}}, \overline{\beta}_i,T) \right> +  \left<\overline{ \beta}_i, \partial_{\beta_i}  f_{ih}(\overline{\boldsymbol{\sigma}},\overline{\beta}_i, T) \right> \Big) +  \lambda \, \partial_T  f_{ih}(\overline{\boldsymbol{\sigma}}, \overline{\beta}_i,T).
\end{aligned}
\end{equation}

\subsubsection{Kinematic hardening} \label{thermo_kin_hard}

For the thermo kinematic hardening law we get the same expressions than in \S \ref{thermo_isot_hard}. The only change is to remplace $\beta_i$ and $\xi_i$ by $\beta_k$ and $\xi_k$, associated with a new constraint $f_{kh}(\boldsymbol{\sigma}, \beta_k,T) \leq 0$ and a new potential function $\mathcal{H}(\xi_k)$ for the kinematic hardening variables.

\subsection{Summary}
 
 We summarise the previous results in the following proposition 
\begin{proposition}
Let the 1D thermoelastoplastic model as described previously. 
Let the total strain $\boldsymbol{\epsilon}$ and the total entropy $S$ which are seen as the sum of their elastic and plastic part, i.e.,
\begin{equation}\label{sum}
\boldsymbol{\epsilon} = \boldsymbol{\epsilon}_e + \boldsymbol{\epsilon}_p, \qquad S=S_e + S_p.
\end{equation}
With the isotropic and kinematic strain hardening variables $\boldsymbol{\xi}=(\xi_i, \xi_k)$, the temperature $T$, and the Helmhotz free energy $\Psi_{h} (\boldsymbol{\epsilon}_e, \boldsymbol{\xi}, T)$. 
Given the yield criterion $f_h(\boldsymbol{\sigma}, \boldsymbol{\beta} ,T) \leq 0$
which constraint the stress tensor field $\boldsymbol{\sigma} = \partial_{\boldsymbol{\epsilon}_e}\Psi_{h}(\boldsymbol{\epsilon}_e, \boldsymbol{\xi}, T)$, the stress hardening variables $\boldsymbol{\beta} =- \partial_{\boldsymbol{\xi}}\Psi_{h}(\boldsymbol{\epsilon}_e, \boldsymbol{\xi}, T)$, and the temperature $T$ to lie in the thermoelastoplastic domain. 
If we assume that there are no heat flux nor heat supply in our model. 
Then, at fixed temperature $T$, when plastic phenomenon occurs at time $\overline{t}$ the elastic and plastic entropy evolve by jumps, and their time rate of change are given by
\begin{equation} \label{time_rate_change}
\begin{aligned}
\left\llbracket S_e \right\rrbracket_t = & \frac{\lambda}{T} \Big( \left< \overline{\boldsymbol{\sigma}}, \partial_{\boldsymbol{\sigma}}  f_{h}(\overline{\boldsymbol{\sigma}}, \overline{\boldsymbol{\beta}},T) \right> +  \left<\overline{ \boldsymbol{\beta}}, \partial_{\boldsymbol{\beta}}  f_{h}(\overline{\boldsymbol{\sigma}},\overline{\boldsymbol{\beta}}, T) \right> \Big) ,
\\
\left\llbracket S_p \right\rrbracket_t  = & \lambda \, \partial_T  f_{h}(\overline{\boldsymbol{\sigma}}, \overline{\boldsymbol{\beta}},T), \quad \text{with} \;\; \lambda> 0.
\end{aligned}
\end{equation}

\end{proposition}

\begin{remark}{\rm
We recall that: ``For every admissible process in a perfect material, the entropy production is zero'' (\cite{TrNo1965}). In our case the entropy production is only due to plastic phenomenon. \quad $\square$ }
\end{remark}

\section{Conclusion}

As highlighted in this development, the elastoplasticity is a nonsmooth phenomenon described by a succession of dissipation jumps which interrupt the smooth path accounted for by a multisymplectic variational formulation. By opposition with viscoelastoplastic dissipative phenomena which are smooth and not described through a variational formulation. 

\medskip

Hence, the next important task is to develop discrete mechanics for nonsmooth elastoplasticity by taking advantage of the variational integrators (such as \cite{FeMaOrWe2003} and \cite{DGBDRA2017}) that are developing in that direction. This task is presently in progress in \cite{Demoures2018c}. 

%In \cite{Demoures2017b} unilateral low impacts are studied through a thermoviscoelastoplastic model.

There are several other directions to pursue. The most important is to include friction in the nonsmooth multisymplectic variational formalism, which is a dissipative phenomenon defined through a maximal principle in the same way as elastoplasticity. 

Then we need to associate different nonsmooth problems, like contact with plasticity, or friction with plasticity, or even contact, plasticity and friction that require further attention in order to get a clear picture of these associations.

\paragraph{Acknowledgment.} I thank Doc. F. Gay-Balmaz for having welcomed me during 6 months at LMD/IPSL, CNRS, Ecole Normale Sup\'erieure Paris.
Moreover I wish to thank Prof. A. Curnier (EPFL) for many helpful discussions.

% \section{Compliance with Ethical Standards}
 
% I attest that :
 
%  \begin{itemize}
%\item
%The manuscript has not been submitted to more than one journal for simultaneous consideration.

%\item The manuscript has not been published previously (partly or in full).

%\item A single study is not split up into several parts to increase the quantity of submissions and submitted to various journals or to one journal over time (e.g. Òsalami-publishingÓ).

%\item No data have been fabricated or manipulated (including images) to support my conclusions.

%\item No data, text, or theories by others are presented as if they were the authorÕs own (ÒplagiarismÓ). 

%\item Disclosure of potential conflicts of interest.

%\end{itemize}
 
{\footnotesize

\bibliographystyle{new}

}

\end{document}